\documentclass{conm-p-l}

\usepackage{amssymb,enumerate}

\usepackage[all]{xy}
\SelectTips{eu}{}

\newcommand{\aq}[4]{\operatorname{D}_{#1}(#2\hskip.7pt|\hskip.7pt#3;#4){}}
\newcommand{\aqc}[4]{\operatorname{D}^{#1}(#2\hskip.7pt|\hskip.7pt#3;#4){}}

\newcommand{\der}[3]{{\operatorname{Der}_{#1}(#2\,;#3)}}
\newcommand{\kahler}[2]{{\Omega_{#1\var #2}}}

\newcommand{\sctan}[1]{{\CL_{#1}}}

\newcommand{\free}[2]{{#1[#2]}}
\newcommand{\add}[2]{{#1[\{#2\}]}}

\newcommand{\koszul}[2]{{\operatorname{K}[#1\,;#2]}}

\newcommand{\norm}[1]{{\mathsf{N}(#1)}}
\newcommand{\simp}[1]{{\mathsf{s}(#1)}}

\newcommand{\sst}{\scriptstyle}
\newcommand{\sss}{\scriptscriptstyle}

\newcommand{\ges}{{\sss\geqslant}}

\newcommand{\bsr}{{\boldsymbol r}}

\newcommand{\eps}{\epsilon}

\newcommand{\card}{\operatorname{card}}

\newcommand{\cls}[1]{{\operatorname{cls}(#1)}}
\newcommand{\col}{\colon}
\newcommand{\dd}{\partial}

\newcommand{\EH}[3]{{}^{#1}\!\operatorname{E}_{#2,#3}}

\newcommand{\fd}{\operatorname{fd}}

\newcommand{\fm}{{\mathfrak m}}
\newcommand{\fn}{{\mathfrak n}}
\newcommand{\fp}{{\mathfrak p}}
\newcommand{\fq}{{\mathfrak q}}

\newcommand{\hH}{\operatorname{H}}
\newcommand{\HH}[2]{\operatorname{H}_{#1}(#2)}
\newcommand{\htpy}[2]{{\pi_{#1}(#2)}}

\newcommand{\hra}{\hookrightarrow}
\newcommand{\id}{{\operatorname{id}}}
\newcommand{\image}{\operatorname{Im}}
\newcommand{\Ker}{\operatorname{Ker}}

\newcommand{\lra}{\longrightarrow}

\newcommand{\pd}{\operatorname{pd}}

\newcommand{\rank}{\operatorname{rank}}

\newcommand{\Spec}{\operatorname{Spec}}

\newcommand{\susp}{{\sst {\mathsf\Sigma}}}

\newcommand{\Ext}[4]{\operatorname{Ext}^{#1}_{#2}(#3,#4){}}
\newcommand{\Hom}[3]{\operatorname{Hom}_{#1}(#2,#3)}

\newcommand{\Tor}[4]{\operatorname{Tor}_{#1}^{#2}(#3,#4){}}

\newcommand{\tra}{\twoheadrightarrow}

\newcommand{\var}{{\hskip.7pt\vert\hskip.7pt}}

\newcommand{\vf}{{\varphi}}

\newcommand{\wh}{\widehat}
\newcommand{\wt}{\widetilde}

\newcommand{\xra}{\xrightarrow}

\newcommand{\BZ}{{\mathbb Z}}

\newcommand{\CL}{{\mathcal L}}

\newcommand{\Mod}{{\mathcal{M}}}

\newcommand{\ST}{{\mathsf T}}

\theoremstyle{plain}

\newtheorem{theorem}{Theorem}[section]
\newtheorem{proposition}[theorem]{Proposition}
\newtheorem{lemma}[theorem]{Lemma}

\newtheorem*{Lemma}{Lemma}
\newtheorem*{Proposition}{Proposition}

\theoremstyle{definition}

\newtheorem{example}[theorem]{Example}

\newtheorem{exercise}[theorem]{Exercise}

\newtheorem{subexercise}{Exercise}[theorem]

\newtheorem{chunk}[theorem]{}

\newenvironment{itchunk}{\begin{chunk}\textit}{\end{chunk}}

\theoremstyle{remark}

\newtheorem{remark}[theorem]{Remark}
\newtheorem{construction}[theorem]{Construction}

\newtheorem{vista}[theorem]{Notes}

\newtheorem*{Question}{Question}

\numberwithin{equation}{theorem}

\begin{document}

\title[Andr\'e-Quillen homology]{Andr\'e-Quillen homology of commutative algebras}

\author{Srikanth Iyengar}
\address{Department of Mathematics,
University of Nebraska, Lincoln, NE 68588, U.S.A.}
\email{iyengar@math.unl.edu}

\thanks{Supported by the National Science Foundation, under grant DMS 0442242} 

\copyrightinfo{2006}{American Mathematical Society}

\date{\today}

\begin{abstract}
  These notes are an introduction to basic properties of Andr\'e-Quillen homology for
  commutative algebras.  They are an expanded version of my lectures at the summer school:
  Interactions between homotopy theory and algebra, University of Chicago, 26th July - 6th
  August, 2004. The aim is to give fairly complete proofs of characterizations of smooth
  homomorphisms and of locally complete intersection homomorphisms in terms of vanishing
  of Andr\'e-Quillen homology. The choice of the material, and the point of view, are
  guided by these goals.
\end{abstract}

\maketitle

\tableofcontents

\section{Introduction}
In the late 60's Andr\'e and Quillen introduced a (co)-homology theory for commutative
algebras that now goes by the name of Andr\'e-Quillen (co)-homology. This is the subject
of these notes. They are no substitute for either the panoramic view that \cite{Qu:ams}
provides, or the detailed exposition in \cite{Qu:mit} and \cite{An:hca}.

My objective is to provide complete proofs of characterizations of two important classes
of homomorphisms of noetherian rings: regular homomorphisms and locally complete
intersection homomorphisms, in terms of Andr\'e-Quillen homology.  However, I have chosen
to treat only the case when the homomorphism is essentially of finite type; this notion is
recalled a few paragraphs below.  One reason for this choice is that it is this class of
homomorphisms which is of principal interest from the point of view of algebraic geometry.

The main reason is that there are technical hurdles, even at the level of definitions and
which have nothing to do with Andr\'e-Quillen homology, that have to be crossed in dealing
with general homomorphisms, and delving into those aspects would be too much of a
digression.  The problem is intrinsic: There are many results for homomorphisms
essentially of finite type (notably, those involving completions and localizations) that
are simply not true in general, and require additional hypotheses. Some of these issues are
discussed in the text.

Andr\'e-Quillen homology is also discussed in Paul Goerss' notes in this volume.  There it
appears as the derived functor of abelianization, while here it viewed as the derived
functor (in a non-abelian context) of K\"ahler differentials.  Another difference is that
in the former, as in Quillen's approach, simplicial resolutions are treated in the general
context of cofibrant replacements in model categories. Here I have described, as Andr\'e
does, an explicit procedure for building simplicial resolutions. This approach is ad
hoc, but it does allow one to construct resolutions in the main cases of interest in these
notes. In any event, it was useful and entertaining to work with `concrete' simplicial
algebras and modules. However, when it comes to establishing the basic properties of
Andr\'e-Quillen homology, I have followed Quillen's more homotopy theoretic treatment, for
I believe that it is easier to grasp.

A few words now on the exposition: Keeping in line with the aim of the summer school, and
the composition of its participants, I have written these notes for an audience of
homotopy theorists and (commutative) algebraists. Consequently, I have taken for granted
material that will be familiar to mathematicians of either persuasion, but have attempted
to treat with some care topics that may be unfamiliar to one or the other. For instance, I
have not hesitated to work with the homotopy category of complexes of modules, and even
its structure as a triangulated category, but I do discuss in detail simplicial
resolutions (presumably for algebraists), and K\"ahler differentials (presumably for
homotopy theorists).

\subsection*{Acknowledgements}
I should like to thank the organizers: Lucho Avramov, Dan Christensen, Bill Dwyer, Mike
Mandell, and Brooke Shipley for giving me an opportunity to speak on Andr\'e-Quillen
homology. I owe special thanks to Lucho Avramov also for numerous discussions on this
writeup.

\subsection*{Notation} The rings in the paper are commutative. 

\begin{itchunk}{Complexes.}
  For these notes, the principal reference for homological algebra of complexes is
  Weibel's book \cite{Wi}, and sometimes also the article \cite{AF}, by Avramov and Foxby.
  Complexes of modules will be graded homologically:
\[
\cdots \lra M_{i+1}\lra M_i \lra M_{i-1}\lra \cdots\,.
\]
When necessary, the differential of a complex $M$ is denoted $\dd^M$.  The
\emph{suspension} of $M$, denoted $\susp M$, is the complex with
\[
(\susp M)_n = M_{n-1} \quad \text{and} \quad \dd^{\susp M} = - \dd^M\,.
\]
Given complexes of $R$-modules $L$ and $M$, the notation $L\simeq M$ indicates that
$L$ and $M$ are homotopy equivalent.
\end{itchunk}

\begin{itchunk}{Homomorphisms.}
  Let $\vf\col R\to S$ be a homomorphism of commutative rings.  One says that $\vf$ is
  \emph{flat} if the $R$-module $S$ is flat.  If the $R$-algebra $S$ is finitely
  generated, then $\vf$ is of \emph{finite type}; it is \emph{essentially of finite type}
  if $S$ is a localization, at a multiplicatively closed set, of a finitely generated
  $R$-algebra.

The notation $(R,\fm,k)$ denotes a (commutative, noetherian) local ring $R$, with maximal
ideal $\fm$, and residue field $k=R/\fm$.  A homomorphism of local rings $\vf\col
(R,\fm,k)\to (S,\fn,l)$ is \emph{local} if $\vf(\fm)\subseteq \fn$. 

For every prime ideal $\fp$ in $R$, we set $k(\fp)=R_\fp/\fp R_\fp$; this is the residue
field of $R$ at $\fp$. The \emph{fiber} of $\vf$ over $\fp$ is the $k(\fp)$-algebra
$S\otimes_Rk(\fp)$.  Given a prime ideal $\fq$ in $S$, the induced local homomorphism
$R_{\fq\cap R}\to S_\fq$ is denoted $\vf_\fq$.

For results in commutative ring theory, we usually refer to Matsumura \cite{Ma}.
\end{itchunk}

\section{K\"ahler differentials}
\label{Kahler differentials}

Let $\vf\col R\to S$ be a homomorphism of commutative rings and $N$ an $S$-module.  

The ring $S$ is commutative, so any $S$-module (be it a left module or a right module) is
canonically an $S$-bimodule; for instance, when $N$ is a left $S$-module, the right
$S$-module structure is defined as follows: for $n\in N$ and $s\in S$, set
\[
n\cdot s = sn 
\]
In what follows, it will be assumed tacitly that any $S$-module, in particular, $N$, is an
$S$-bimodule, and hence also an $R$-bimodule, via $\vf$.

\begin{itchunk}{Derivations.}
\label{derivations}  
An \emph{$R$-linear derivation} of $S$ with coefficients in $N$ is a homomorphism of $R$-modules
$\delta\col S\to N$ satisfying the Leibniz rule: 
\[
 \delta(st) = \delta(s)t + s\delta(t) \quad\text{for $s,t\in S$.}
\]
An alternative definition is that $\delta$ is a homomorphism of abelian groups satisfying
the Leibniz rule, and such that $\delta\vf=0$.  The set of $R$-linear derivations of $S$
with coefficients in $N$ is denoted $\der RSN$.  This is a subset of $\Hom RSN$, and even
an $S$-submodule, with the induced action:
\[
(s\cdot \delta)(t) =s\delta(t) 
\]
for $s,t\in S$ and $\delta\in \der RSN$. 

\begin{subexercise}
\label{derivations:hom}
  Let $M$ be an $S$-module. The  homomorphism of $S$-modules
\begin{gather*}
\Hom SMN\otimes_S \Hom RSM \lra \Hom RSN\\
\alpha\otimes \beta \mapsto \alpha\beta 
\end{gather*}
restricts to a homomorphism of $S$-modules:
\[
\Hom SMN \otimes_S \der RSM \lra \der RSN
\]
In particular, for each derivation $\delta\col S\to M$, composition induces a homomorphism
of $S$-modules $\Hom SMN \to \der RSN$.
\end{subexercise}
\end{itchunk}

\begin{itchunk}{K\"ahler differentials.}
\label{kahler:construction}
It is not hard to verify that the map $N\mapsto \der RSN$ is an additive functor on
the category of $S$-modules.  It turns out that this functor is representable, that is to
say, there is an $S$-module $\Omega$ and an $R$-linear derivation $\delta\col S\to \Omega$
such that, for each $S$-module $N$, the induced homomorphism 
\begin{gather*}
\label{der:representable}
\Hom S{\Omega}N \lra \der RSN
\end{gather*}
of $S$-modules, is bijective. Such a pair $(\Omega,\delta)$ is unique up to isomorphism, in
a suitable sense of the word; one calls $\Omega$ the module of \emph{K\"ahler
  differentials} and $\delta$ the \emph{universal derivation} of $\vf$. In these note,
they are denoted $\Omega_\vf$ and $\delta^\vf$ respectively; we sometimes follow
established usage of writing $\kahler SR$ for $\Omega_\vf$.

In one case, the existence of such an $\Omega_\vf$ is clear:

\begin{subexercise}
Prove that when $\vf$ is surjective $\Omega_\vf=0$.
\end{subexercise}

A homomorphism $\vf$ such that $\Omega_\vf=0$ is said to be \emph{unramified}.

There are various constructions of the module of K\"ahler differentials and the universal
derivation; see Matsumura \cite[\S9]{Ma}, and Exercise \ref{kahler:presentation}.
The one presented below is better tailored to our needs:

We are in the world of commutative rings, so the product map
\[
\mu_R^S\col S\otimes_RS\lra S \quad\text{where $s\otimes t\mapsto st$}
\]
is a homomorphism of rings. Set $I=\Ker(\mu_R^S)$. Via $\mu_R^S$ the $S\otimes_RS$-module
$I/I^2$ acquires the structure of an $S$-module. Set
\begin{gather*}
\Omega_\vf = I/I^2 \quad \text{and}\quad 
\delta^\vf \col S\lra \Omega_\vf \quad\text{with $\delta^\vf(s)= (1\otimes s - s\otimes 1)$.}
\end{gather*}
As the notation suggests, $(\Omega_\vf,\delta^\vf)$ is the universal pair we seek. The
first step in the verification of this claim is left as an

\begin{subexercise}
The map  $\delta^\vf$ is an $R$-linear derivation.
\end{subexercise}

By Exercise \ref{derivations:hom}, the map $\delta^\vf$ induces a homomorphism of
$S$-modules
\begin{equation*}
 \Hom S{\Omega_\vf}N \lra \der RSN  \tag{$\ast$}
\end{equation*}

We prove that this map is bijective by constructing an explicit inverse.

Let $\delta\col S\to N$ be an $R$-linear derivation. As $\delta$ is a homomorphism of $R$-modules,
extension of scalars yields a homomorphism of $S$-modules
\[
\delta'\col S\otimes_RS\to N, \quad\text{where}\quad \delta'(s\otimes t) = s\delta(t) 
\] 
Here we view $S\otimes_RS$ as an $S$-module via action on the left hand factor of the
tensor product: $s\cdot (x\otimes y)=(sx\otimes y)$. One thus obtains, by
restriction, a homomorphism of $S$-modules $I\to N$, also denoted $\delta'$.

\begin{subexercise} Verify the following claims.
\begin{enumerate}[\quad\rm(1)]
\item
$\delta'(I^2)=0$, so $\delta'$ induces a homomorphism of $S$-modules
  $\wt\delta \col\Omega_\vf \to N$.
\item The assignment $\delta\mapsto \wt\delta$ gives a homomorphism $\der RSN\to \Hom
  S{\Omega_\vf}N$ of $S$-modules, and it is an inverse to the map ($\ast$) above.
\end{enumerate}
\end{subexercise}

This exercise justifies the claim that $\delta^\vf\col S\to \Omega_\vf$ is a universal
derivation.
\end{itchunk}

The next goal is an explicit presentation for $\Omega_\vf$ as an $S$-module, given the
presentation of $S$ as an $R$-algebra. The first step towards it is the following
exercise describing the module of K\"ahler differentials for polynomial extensions of $R$.
Solve it in two ways: by using the construction in paragraph \ref{kahler:construction}
above; by proving directly that it has the desired universal property.

\begin{exercise}
\label{kahler:freeextension}
Let $S=R[Y]$ be the polynomial algebra over $R$ on a set of variables $Y$, and 
$\vf\col R\to S$ the inclusion map.  Prove that
\begin{gather*}
\Omega_\vf = \bigoplus_{y \in Y}S{dy} \quad \text{and}\quad
\delta^{\vf}(r) = \sum_{y\in Y} \frac{\dd{r}}{\dd{y}}dy\,.
\end{gather*}
Here, $\Omega_\vf$ is a free $S$-module on a basis $\{dy\}_{y\in Y}$, and
$\dd(-)/\dd{y}$ denotes partial derivative with respect to $y$.
\end{exercise}

\begin{itchunk}{Jacobi-Zariski sequence.}
  Let $Q\xra{\psi} R\xra{\vf}S$ be a homomorphism of commutative rings.  One has a natural
  exact sequence of $S$-modules:
\begin{equation}
\label{jz:kahler}
\xymatrixcolsep{2pc}
\xymatrix{
\Omega_{\psi}\otimes_RS  \ar@{->}[r]^-{\alpha}
  &\Omega_{\vf\psi}\ar@{->}[r]^-{\beta}
  &\Omega_{\vf} \ar@{->}[r]& 0\,.}
\end{equation}
The maps in question are defined are follows: by restriction, the $R$-linear derivation
$\delta^\vf\col S\to \Omega_\vf$ is also a $Q$-linear derivation, hence it induces the
homomorphism $\beta\col \Omega_{\vf\psi}\to \Omega_{\vf}$ such that
$\beta\circ\delta^{\vf\psi}=\delta^\vf$. In the same vein, $\delta^{\vf\psi}\vf\col R\to
\Omega_{\vf\psi}$ is a $Q$-linear derivation, so it induces an $R$-linear homomorphism
$\alpha'\col \Omega_{\psi}\to \Omega_{\vf\psi}$; the map $\alpha$ is obtained by extension
of scalars, for $\Omega_{\vf\psi}$ is an $S$-module.

I leave it to you to verify that the sequence \eqref{jz:kahler} is exact. It is sometimes
called the Jacobi-Zariski sequence. One way to view Andr\'e-Quillen homology is that it
extends this exact sequence to a long exact sequence; that is to say, it is a `left
derived functor' of $\Omega_{-}$, viewed as a functor of algebras; see
\ref{transitivity}.

When $N$ is a $S$-module, applying $\Hom S-N$ to the exact sequence \eqref{jz:kahler}, and
using the identification in \ref{der:representable}, yields an exact sequence of
$S$-modules
\[
0\lra \der RSN \lra  \der QSN \lra  \der QRN \,.
\]
One could just as well have deduced \eqref{jz:kahler} from the naturality of this sequence.

\begin{subexercise}
Interpret the maps in the exact sequence above.    
\end{subexercise}
\end{itchunk}

The following exercise builds on Exercise \ref{kahler:freeextension}.

\begin{exercise}
\label{kahler:freeextensionmaps}
 Let $\psi\col R[Y]\to R[Z]$ be a homomorphism of $R$-algebras, where $Y$ and $Z$ are
sets of variables. Verify that the map 
\begin{gather*}
\kahler {\psi}R\col \kahler {R[Y]}R \lra \kahler {R[Z]}R\quad \text{is defined by} \quad
   y \mapsto \sum_{z\in Z} \frac{\dd{y}}{\dd{z}}dz\,.
\end{gather*}

\end{exercise}

Sequence \eqref{jz:kahler} allows for a `concrete' description of the K\"ahler
differentials:

\begin{exercise}
\label{kahler:presentation}
Write $S=R[Y]/(\bsr)$, where $R[Y]$ is the polynomial ring over $R$ on a set of variables
$Y$ and $\bsr =\{r_\lambda\}$ is a set of polynomials in $R[Y]$ indexed by the set
$\Lambda$. Note that such a presentation of $S$ is always possible.  The homomorphism
$\vf$ is then the composition $R\hra R[Y]\tra S$.

Prove that the module of K\"ahler differentials of $\vf$ is presented by
\begin{align*}
  \bigoplus_{\lambda\in\Lambda} Se_\lambda &\xra{\dd} 
 \bigoplus_{y\in Y} Sdy \lra \Omega_\vf\lra 0\,,
  \quad\text{where}\\
  &\dd(e_\lambda) = \sum_{y\in Y} \frac{\dd(r_\lambda)}{\dd{y}}dy\,.
\end{align*}
The matrix representing $\dd$ is the Jacobian matrix of the polynomials $\bsr$.
\end{exercise}

Here is an exercise to give you a feel for the procedure outlined above:

\begin{exercise}
  Let $k$ be a field, $S=k[y]$, the polynomial ring in the variable $y$, and let $R$ be
  the subring $k[y^2,y^3]$. Find a presentation for the module of K\"ahler differentials
  for the inclusion $R\hra S$.
\end{exercise}

The relevance of the following exercise should be obvious.

\begin{exercise}
\label{kahler:products}
Let $S$ and $T$ be $R$-algebras. Prove that there is a natural homomorphism of 
$(S\otimes_RT)$-modules
\[
(S\otimes_R\kahler TR) \bigoplus (\kahler SR \otimes_R T) \lra
\kahler {(S\otimes_RT)}R\,,
\]  
and that this map is bijective.
\end{exercise}

In a special case, one can readily extend \eqref{jz:kahler} one step further to the left:

\begin{itchunk}{The conormal sequence.}
  Suppose $S=R/I$, where $I$ is an ideal in $R$, and $\vf\col R\to S$ the canonical
  surjection; in particular, $\Omega_{\vf}=0$.  Let $\psi\col Q\to R$ be a homomorphism of
  rings. The exact sequence \eqref{jz:kahler} extends to an exact sequence 
\[
\xymatrixcolsep{2pc}
\xymatrix{
I/I^2\ar@{->}[r]^-{\zeta} & \Omega_{\psi}\otimes_RS  \ar@{->}[r]^-{\alpha}
  &\Omega_{\vf\psi}\ar@{->}[r] & 0}
\]
of $S$-modules.  The map $\zeta$ is defined as follows: restricting the universal
derivation $\delta^\psi$ gives a $Q$-linear derivation $I\to \Omega_{\psi}$ and hence, by
composition, a $Q$-linear derivation $\delta\col
I\to(\Omega_{\psi}\otimes_RS)=\Omega_{\psi}/I\Omega_{\psi}$.  Keeping in mind that
$\delta$ is a derivation it is easy to verify that $\delta(I^2)=0$, so it factors through
$I/I^2$; this is the map $\zeta$.  It is also elementary to check that $\zeta$ is
$S$-linear.
\end{itchunk}

\section{Simplicial algebras}
\label{Simplicial algebras}
This section is a short recap on simplicial algebras and simplicial modules.  The aim is
to introduce enough structure, terminology, and notation to be able to work with
simplicial algebras and their resolutions, and construct cotangent complexes, the topics
of forthcoming sections. The reader may refer to \cite{GJ} and \cite{May} for in-depth
treatments of things simplicial.

To begin with, let me try to explain what we are trying to do here.

\begin{itchunk}{Computing derived functors.}
  I remarked during the discussion on Jacobi-Zariski sequence \eqref{jz:kahler} that
  Andr\'e-Quillen homology may be viewed as a derived functor of $\Omega_{-}$. In order to
  understand the problem, and its solution, let us revisit the process of deriving a more
  familiar functor.

  As before, let $\vf\col R\to S$ be a homomorphism of commutative rings.  Let $\Mod(S)$
  be the category of $S$-modules. Consider the functor
\[
S\otimes_R -\col \Mod(S)\lra \Mod(S) \quad \text{where $N\mapsto S\otimes_RN$.}
\]
Each exact sequence of $S$-modules $0\to N\to N'\to N''\to 0$ gives rise to an exact
sequence of $S$-modules
\[
\tag{$\ast$}
S\otimes_RN\lra S\otimes_RN'\lra S\otimes_RN''\lra 0
\] 
However, the homomorphism on the left is not injective, unless $S$ is flat as an
$R$-module; in short, the functor $S\otimes_R-$ is left-exact, but it is not exact.  This
lack of exactness is compensated by extending the sequence above to a long exact sequence.
There are three steps involved in this process:

\smallskip

\emph{Step }1. Construct a projective resolution $F$ of $S$ over $R$.

\smallskip

\emph{Step} 2. Show that $F$ is unique up to homotopy of complexes of $R$-modules.

\smallskip

\emph{Step} 3. Set $\ST_\vf= F\otimes_RS$, and for each $S$-module $N$ set
\[
\hH^\vf_n(N) = \HH n{\ST_\vf\otimes_SN}
\]
One then has $\hH^\vf_0(N)= S\otimes_RN$, and the functors $\{\hH^{\vf}_n(-)\}_{n\ges 0}$
extend ($\ast$) above to a long exact sequence of $S$-modules, which is what one wants.
We have not discovered anything new here: $\hH^\vf_n(N) = \Tor nRSN$.

Note that the complex $\ST_\vf$ is well-defined in the homotopy category of complexes of
$S$-modules; this follows from Step 2. It is this property that dictates the kind of
resolution we pick. For instance, flat resolutions, although a natural choice, would not
work, for they are not unique, even up to homotopy.  It is another matter that they can be
used to compute $\hH^\vf_n(N)$.

We turn now to the functor of interest $\Omega_{-}$. Taking a cue from the preceding
discussion, the plan is to attach a complex of projective $S$-modules called the
\emph{cotangent complex of $\vf$}, which I denote $\sctan \vf$, with
\[
\HH 0{\sctan \vf \otimes_SN} = \Omega_{\vf}\otimes_SN
\]
such that it extends the sequence \eqref{jz:kahler} to a long exact sequence of
$S$-modules. The functor $\Omega_{-}$ is non-linear: it takes into account the structure
of $S$ as an $R$-algebra, rather than as an $R$-module. Keeping this in mind, one should
pick a suitable category of $R$-algebras, and a notion of homotopy for morphisms in that
category, such that resolutions have the following properties:
\begin{enumerate}[\quad\rm(a)]
\item They must reflect the structure of $S$ as an $R$ algebra.
\item They should be unique up to homotopy.
\item The functor $\Omega_{-}$ must preserve homotopies, in a suitable sense of the word.
\end{enumerate}

It turns out that simplicial algebras provide the right context for obtaining such resolutions;
confer \cite{Pg:uc} for a discussion about why this is so. 

\begin{Question}
Why is the category of differential graded $R$-algebras not suitable for the purpose on hand?  
\end{Question}
  
\end{itchunk}

\begin{itchunk}{Simplicial modules and algebras.}
\label{Simplicial objects}
As usual, let $R$ be a commutative ring.

A \emph{simplicial $R$-module} is a simplicial object in the category of $R$-modules, that
is to say, a collection $V=\{V_n\}_{n\ges 0}$ of $R$-modules such that, for each
non-negative integer $n$, there are homomorphisms of $R$-modules:
\[
d_i\col V_n\to V_{n-1} \quad \text{and} \quad s_j\col V_n\to V_{n+1} 
\qquad\text{for $0\leq i,j\leq n.$}
\]
called \emph{face maps} and \emph{degeneracies}, respectively, satisfying the identities:
\begin{align}
\label{simplicial:identities}
&d_id_j = d_{j-1}d_i  \quad \text{when $i<j$} \\ \notag
&d_is_j = \begin{cases}
s_{j-1}d_i &  \text{when $i<j$} \\
1          & \text{when $i=j,j+1$} \\
s_jd_{i-1} & \text{when $i>j+1$}
           \end{cases} \\ \notag
&s_is_j = s_{j+1}s_i \quad \text{when $i\leq j$} \notag
\end{align}
Prescribing this data is equivalent to defining a contravariant functor from the ordinal
number category to the category of $R$-modules; see \cite[\S(1.5)]{Pg:uc}

A \emph{simplicial $R$-algebra} is a simplicial object in the category of $R$-algebras;
thus, it is a simplicial $R$-module $A$ where each $A_n$ has the structure of an
$R$-algebra, and the face and degeneracies are homomorphisms of $R$-algebras.  A
\emph{simplicial module} over a simplicial $R$-algebra $A$ is a simplicial $R$-module $V$ where
each $V_n$ is an $A_n$-module and the face maps and degeneracies are compatible
with those on $A$.
\end{itchunk}

\begin{example}
\label{simplicial:rings}
Given an $R$-module $N$, it is not hard to verify that the graded module $\simp N$, with
\[
\simp N_n= N \quad \text{and}\quad d_i = \id^N = s_j \qquad \text{for $0\leq i,j\leq n$.}
\]
is a simplicial $R$-module. For any $R$-algebra $S$, it is evident that $\simp S$ is a
simplicial $R$-algebra.
\end{example}

\begin{itchunk}{Normalization}.  
\label{moore}
Let $V$ be a simplicial $R$-module.  The \emph{normalization} of $V$ is the complex of
$R$-modules $\norm V$, defined by
\begin{align*}
&\norm V_n = \bigcap_{i=1}^{n}\Ker(d_i) \quad\text{with differential}\\
&\dd_n=d_0 \col \norm V_n \lra \norm V_{n-1}\,.
\end{align*}
That $\dd_n$ is a differential follows from \eqref{simplicial:identities}.  The
\emph{$n$th homotopy module} of $V$ is the $R$-module
\[
\htpy nV = \HH nV\,.
\]
This is not the `right' way to introduce homotopy, but will serve the purpose here, see
\cite[(2.15)]{Pg:uc}.

There is another way to pass from simplicial modules to complexes: The face maps on $V$
give the graded $R$-module underlying $V$ the structure of a complex of $R$-modules, with
differential:
\[
\dd_n = \sum_{i=0}^n (-1)^i d_i \col V_n \to V_{n-1}\,.
\]
This complex is also denoted $V$; this could cause confusion, but will not, for the
structure involved is usually clear from the context.  Fortunately, the homology of this
complex is the same as the homotopy, see, for instance, \cite[Chapter III, (2.7)]{GJ}.
\end{itchunk}

\begin{exercise}
\label{simp:norm}
Let $S$ be an $R$-algebra.  Verify that $\norm{\simp S}=S$.
\end{exercise}

\begin{exercise}
\label{pi:properties}
 Let $A$ be a simplicial $R$-algebra and $V$ a simplicial $A$-module. 

 Verify that for each $n$, the $R$-submodules $\norm V_n$, $\Ker(\dd_n)$, and
 $\image(\dd_{n+1})$ of $V_n$ are stable under the action of $A_n$, that is to say, they
 are $A_n$-submodules of $V_n$.

 In particular, $\htpy 0A$ is an $R$-algebra. Moreover, when $V_n$ is a noetherian
 $A_n$-module, so is $\htpy nV$.
\end{exercise}

\begin{vista}
  For each simplicial $R$-algebra $A$, the graded module $\htpy *A$ is a commutative
  $\htpy 0A$-algebra with divided powers, see \cite[(2.3)]{Pg:smf}.  Moreover, if $V$ is a
  simplicial $A$-module, $\htpy *V$ is a graded $\htpy *A$-module.  These structures play
  no role in this write-up, but they are an important facet of the theory; see
  \cite{Av:am} and \cite{AI:ens}.
\end{vista}

\begin{vista}
  The functor $\norm{-}$ from simplicial $R$-modules to complexes of $R$-modules is an
  equivalence of categories; this is the content of the Dold-Kan theorem, see
  \cite[(4.1)]{Pg:uc}.
\end{vista}

\begin{example}
\label{modules:smodules}
Let $A$ be a simplicial $R$-algebra and $N$ an $\htpy 0A$-module. Then $\simp N$ is a
simplicial $A$-module, where, for each non-negative integer $n$, the $A_n$-module
structure on $\simp N_n$ is induced via the composed homomorphism of rings
\[
A_n\xra{d_{i_n}}A_{n-1}\xra{d_{i_{n-1}}}\cdots \xra{d_{i_1}} A_0\,.
\]
It is an exercise to check that the choices of indices $i_n,\dots,i_1$ is irrelevant.
\end{example}

\begin{itchunk}{Morphisms.}
A \emph{morphism} $\Phi\col A\to B$ of simplicial $R$-algebras is a collection of
homomorphisms of $R$-algebras $\Phi_n\col A_n\to B_n$, one for each $n\ge 0$, commuting
with both face maps and degeneracies. Such a $\Phi$ induces a homomorphism of $R$-modules
\[
\htpy {*}\Phi \col \htpy {*}A \to \htpy *B
\]
One says that $\Phi$ is a \emph{weak equivalence} if $\htpy *\Phi$ is bijective.  I will
leave it to you to formulate the definition of a morphism of simplicial modules.
\end{itchunk}

\begin{example}
\label{morphism:rings}
Let $A$ be a simplicial $R$-algebra.  Following Example \ref{modules:smodules}, it is not
hard to verify that any homomorphism of $R$-algebras $\phi\col \htpy 0A\to S$ induces a
morphism of simplicial algebras, $\Phi \col A\to \simp S$. Given Example
\ref{simplicial:rings} it is clear that
\[
\htpy n{\Phi}=\begin{cases}
\phi & \text{when $n=0$}\\
0    & \text{otherwise}
\end{cases}
\]
Thus, $\Phi$ is a weak equivalence if and only if $\phi$ is bijective and $\htpy
nA=0$ for $n\ge 1$.
\end{example}

To summarize Examples \ref{simplicial:rings} and \ref{morphism:rings}: The functor
$\simp{-}$ is a faithful embedding of the category of $R$-algebras into the category
of simplicial $R$-algebras, and $\htpy 0{-}$ is a left adjoint to this embedding.

\begin{subexercise}
  Prove that the embedding is also full.
\end{subexercise}

\begin{itchunk}{Tensor products.}
\label{tensors}
The \emph{tensor product} of simplicial $A$-modules $V$ and $W$ is the simplicial
$A$-module denoted $V\otimes_AW$, with
\[
(V\otimes_AW)_n = V_n\otimes_{A_n} W_n \quad \text{for each $n\geq 0$,}
\]
and face maps and degeneracies induced from those on $V$ and $W$. When $N$ is an
$\htpy 0A$-module, it is customary to write $V\otimes_AN$ for $V\otimes_A\simp N$.

Various standard properties of tensor products (for example: associativity and commutativity)
carry over to this context.
\end{itchunk}

\section{Simplicial resolutions}
\label{Simplicial resolutions}
This section discusses simplicial resolutions.  The first step is to introduce free
extensions, which are analogues in simplicial algebra of bounded-below complexes of free
modules in the homological algebra of complexes over rings.

\begin{itchunk}{Free simplicial extensions.}
  Let $A$ be a simplicial $R$-algebra. We call a \emph{free}\footnote{See footnote for
Definition 4.20 of \cite{Pg:uc}} simplicial extension of $A$ on a graded set $X=\{X_n\}_{n\ges0}$
  of indeterminates a simplicial $R$-algebra, denoted $\free AX$, satisfying the following
  conditions:
  \begin{enumerate}[\quad\rm(i)]
  \item $\free AX_n = A_n[X_n]$, the polynomial ring over $A_n$ on the variables $X_n$;
  \item $s_j(X_n) \subseteq X_{n+1}$ for each $j,n$;
  \item The inclusion $A\hra \free AX$ is a morphism of simplicial $R$-algebras. 
  \end{enumerate}
Note that there is no restriction on the face maps.

  For instance, if $S=R[Y]$ is a polynomial ring over $R$, then the simplicial algebra
  $\simp{S}$ is a free extension of $\simp R$, with $X_n=Y$ for each $n$.
\end{itchunk}

\begin{itchunk}{Base change.}
  If $\free AX$ is a free extension of $A$, and $\Phi\col A\to B$ is a morphism of
  simplicial algebras, then $B\otimes_A\free AX$ is a free extension of $B$.
\end{itchunk}

\begin{itchunk}{Existence of resolutions.}
\label{resolutions:existence}
  Let $\phi\col A\to B$ be a morphism of simplicial $R$-algebras.  A \emph{simplicial resolution}
  of $B$ over $A$ is a factorization of $\phi$ as a diagram
\[
\xymatrixcolsep{2pc}
\xymatrix{
A \ar@{>->}[r] & \free AX \ar@{->>}[r]^{\Phi} & B}
\]
of morphisms of simplicial algebra, with $A\to \free AX$ a free extension and $\Phi$
 a surjective weak equivalence.  Usually, one refers to $\free AX$ itself as a
simplicial resolution of $B$ over $A$. Simplicial resolutions exist; one procedure for
constructing them is described in the paragraphs below; see \ref{resolutionsexist}.
\end{itchunk}

\begin{itchunk}{Existence of lifting.}
\label{free:lifting1}
Given a commutative diagram of simplicial algebras 
\[
\xymatrixrowsep{2pc}
\xymatrixcolsep{2pc}
\xymatrix{
A \ar@{->}[d]\ar@{->}[r]  & B\ar@{->}[d]^{\Phi}_{\simeq} \\
\free AX \ar@{-->}[ur]^{\kappa} \ar@{->}[r]  & C}
\]
where $\Phi$ is surjective and a weak equivalence, there exists a morphism $\kappa$
that preserves the commutativity of the diagram. 
\end{itchunk}

For a proof of the lifting property, see \cite[(5.4)]{Pg:uc}.  The morphism $\kappa$ in
the diagram above is unique up to homotopy, in a sense described below. The definition may
appear to come out of the blue, but it is a special case of a notion of homotopy
in model categories. Much of this following discussion is best viewed in that general
context; see Dwyer and Spalinski \cite[(4.1)]{DS}, or \cite[(2.2)]{Pg:uc}.

\begin{itchunk}{Homotopy.}
\label{homotopy}
Let $A\to \free AX$ be a free extension. For each integer $n$, one has the product morphism 
\[
\mu_n = A_n[X_n]\otimes_{A_n}A_n[X_n]\lra A_n[X_n]\,.
\]
They form a morphism of simplicial $A$-algebras 
\[
\mu\col \free AX\otimes_A\free AX\to \free AX\,.
\]
It is convenient to write $\free A{X,X}$ for $\free AX\otimes_A\free AX$.

The simplicial algebra $\free A{X,X}$ has the functorial properties one expects of a product.
Namely, given morphisms $\Phi,\Psi\col \free AX\to B$ of simplicial $A$-algebras, there is an
induced morphism of simplicial algebras
\[
\Phi\odot \Psi \col \free A{X,X}\lra B\quad\text{with}\quad
(\Phi\odot\Psi)_n(x\otimes x') = \Phi(x)\Psi(x')\,.
\]

Let $\free A{X,X,Y}$ be a simplicial resolution of $\free AX$ over $\free A{X,X}$; it is
called a \emph{cylinder object} for the $A$-algebra $\free AX$, see \cite[(2.4)]{Pg:uc}.
The morphisms $\Phi$ and $\Psi$ are \emph{homotopic} if $\Phi\odot \Psi$ extends to a
cylinder object, that is to say, there is a commutative diagram of morphisms of simplicial
$A$-algebras
\[
\xymatrixrowsep{2.5pc}
\xymatrixcolsep{2.5pc}
\xymatrix{
\free A{X,X} \ar@{->}[dr]_{\Phi\odot\Psi}\ar@{->}[r] & \free A{X,X,Y} \ar@{->}[d] \\
  &B}
\]
Given the lifting property of free extensions \ref{free:lifting1}, it is easy to check
that the notion of homotopy does not depend on the choice of a cylinder object, and that
homotopy is an equivalence relation on morphisms of $A$-algebras, see \cite[(4.7)]{DS}.
\end{itchunk}

\begin{itchunk}{Uniqueness of lifting.}
\label{free:lifting2}
Using the lifting property \ref{free:lifting1} of free extensions a formal argument
shows that the lifting map $\kappa$ in \ref{free:lifting1} is unique up to homotopy of
simplicial $A$-algebras, see \cite[(4.9)]{DS}.
\end{itchunk}

\begin{itchunk}{Uniqueness of resolutions.}
\label{resolutions:uniqueness}
  A standard argument using lifting properties of free extensions \ref{free:lifting1},
  \ref{free:lifting2} yields that simplicial resolutions are unique up to homotopy of
  simplicial $A$-algebras.

Given an $R$-algebra $S$, it is accepted usage to speak of a simplicial resolution of the
$R$-algebra $S$, meaning a simplicial resolution of $\simp S$ over $\simp R$.
\end{itchunk}

\begin{remark}
\label{resolution:underlying}
Let $S$ be an $R$-algebra and $M$ an $R$-module.  Let $\free RX$ be a simplicial
resolution of the $R$-algebra $S$.  The complex underlying $\free RX$ is an $R$-free
resolution of $S$, so for each integer $n$, one has
\[
\htpy n{\free RX\otimes_RM} = \Tor nRSM\,.
\]
\end{remark}

Next we outline a procedure for constructing simplicial resolutions. The strategy is the
one used to obtain free resolutions:

\begin{itchunk}{Resolutions of modules over rings.}
\label{linear resolutions}
Let $M$ be an $R$-module. A free resolution of $M$ over $R$ may be built as follows: One
constructs a sequence of complexes of free $R$-modules $0\subset F^{(0)}\subseteq
F^{(1)}\subseteq \cdots$ such that $F^{(0)}$ is a free module mapping onto $M$, and for
each $d\geq 1$ one has
\[
\HH i{F^{(d)}} = \begin{cases}
M & \text{for $i=0$\,,} \\
0 & \text{for $1\leq i\leq d-1$}\,.
\end{cases}
\]
Given $F^{(d-1)}$ one builds $F^{(d)}$ by killing cycles in $F^{(d-1)}_{d-1}$ that are not
boundaries. In detail: choose a set of cycles $\{z_\lambda\}_{\lambda\in\Lambda}$ which
generate $\HH {d-1}{F^{(d-1)}}$, and set
\begin{gather*}
  F_{d} = \bigoplus_{\lambda\in\Lambda}Re_\lambda\\
  F^{(d)}= F^{(d-1)}\bigoplus \susp^{d}F_{d} \quad \text{with}\quad 
\dd(e_\lambda) = z_\lambda\,.
\end{gather*}
The homology of $F^{(d)}$ is readily computed from the short exact sequence of complexes
$0\to F^{(d-1)}\to F^{(d)}\to \susp^{d}F_{d}\to 0$.

Then the complex $\cup_{d\ges 0}F^{(d)}$ is the desired free resolution of $M$.
\end{itchunk}

A similar procedure can be used to construct simplicial resolutions.  The crucial step
then is a method for killing cycles. In the category of modules, to kill a cycle in degree
$d-1$ we attached a free module, $F_{d}$, in degree $d$. In the category of simplicial
algebras, we has to attach a (polynomial) variable in degree $d$; however, the simplicial
identities \eqref{simplicial:identities} (notably, $d_is_j=1$ for $i=j,j+1$) force us to
then attach a whole slew of variables in higher degrees.

\begin{itchunk}{Killing cycles.}
\label{killingcycles}
Let $A$ be a simplicial $R$-algebra, $d$ a positive integer, and let $w\in A_{d-1}$ be
a cycle in $\norm A_{d-1}$, the normalized chain complex of $A$, see \ref{moore}.

The goal is to construct a free extension of $A$ in which the cycle $w$ becomes a
boundary; I write $\free A{\{x\}\mid \dd(x)=w}$, or just $\add Ax$ when the cycle being
killed is understood, for the resulting simplicial algebra. It has the following
properties:
\begin{enumerate}[\quad\rm(a)]
\item For each integer $n$, the $A_n$-algebra $\add Ax_n$ is free on a set $X_n$ of finite
  cardinality. In particular, if the ring $A_n$ is noetherian, so is the ring $\add Ax_n$.
\item The inclusion $A\hra \add Ax$ induces isomorphisms
\[
\htpy nA\cong \htpy n{\add Ax}\qquad \text{for $n\leq d-2$.}
\]
\item One has an exact sequence of $\htpy 0A$-modules
\[
0\to A_{d-1}\cls w \to \htpy {d-1}A \to \htpy {d-1}{\add Ax}\to 0\,.
\]
where $\cls w$ is the class of the cycle $w$ in $\htpy {d-1}A$.  Note that the ideal
$A_{d-1}w\subseteq A_{d-1}$ consists of cycles.
\end{enumerate}

The construction of $\add Ax$ is as follows: 

\medskip

\textit{The set $X$}.  For each positive integer $n$, set
\[
X_n = \{x_t\mid \text{$t\col [n]\to [d]$ is surjective and monotone.}\}
\]
Clearly, $\card(X_n)$ is finite, as claimed.

The face and degeneracies on $\add Ax$ extend those on $A$, so to define them, it suffices
to specify their action the set $X$. This process involves the co-face and co-degeneracy
maps, see \cite[(1.10)]{Pg:uc}.

\textit{Degeneracies}. For each $x_t\in X_n$, set
\[
s_j(x_t) = x_{t\circ s^j} \quad \text{for $0\leq i\leq n$.}
\]
Here $s^j\col [n+1]\to [n]$ is the $j$th co-degeneracy operator.

\medskip 

\textit{Face maps}.  The set $X_d$ is a singleton: $\{x_\id\}$. Set
\[
d_i(x_\id) = \begin{cases} 
w & i=0 \\
0 & 1\leq i\leq d
\end{cases}
\]
It remains to define face maps on $X_n$ for $n\ge d+1$.  Fix such an $n$ and a surjective
monotone map $t\col [n]\to [d]$. If for a co-face map $d^i\col [n-1]\to [n]$ the composed
map $t\circ d^i \col [n-1]\to [d]$ is not surjective, then one has a commutative diagram
\begin{equation*}
\xymatrixrowsep{2.5pc}
\xymatrixcolsep{2.5pc}
\xymatrix{
  [n] \ar@{->}[r]^t         & [d] \ar@{<-}[d]^{d^{j}} & \\
[n-1] \ar@{->}[u]^{d^i}  & [d-1] \ar@{<-}[l]^{u}}
\end{equation*}
where $u$ is surjective and monotone.  The face operators on $X_n$ is defined as follows:
\[
d_i(x_t) = \begin{cases}
x_{t\circ d^i} & \text{if $t\circ d^i$ is surjective,}\\
u_*(d_{j}(w)) & \text{otherwise.}
\end{cases}
\]
where $u_*\col A_{d-1}\to A_{n-1}$ is the map induced by $u$. Here we are viewing $A$ as
defining a contravariant functor from the ordinal number category to the category of $R$-algebras;
 $u_*$ is the image of $u$ under this functor, see \cite[(1.10)]{Pg:uc}.

\begin{subexercise}
  Prove that, with the prescription above, $\add Ax$ is a simplicial $R$-algebra, and a
  free extension of $A$.
\end{subexercise}

As to the properties of $\add Ax$: given \ref{killingcycles}(a), it is clear that
\ref{killingcycles}(b) holds.

\begin{subexercise}
Prove that property \ref{killingcycles}.(c) holds. Hint: use \ref{moore}.
\end{subexercise}
\end{itchunk}

\begin{itchunk}{Generating cycles.}
\label{generatingcycles}
The preceding construction can also be used to generate cycles in degree $d$: take $w=0$.
\end{itchunk}

Using \ref{killingcycles} and \ref{generatingcycles}, and taking limits one arrives at the
conclusion below.  This is a good place to point out that this construction of
resolutions is best viewed in the context of a general technique called the `small object
argument', due to Quillen, see \cite[(3.5)]{Pg:uc}.

\begin{itchunk}{Resolutions exist.}
\label{resolutionsexist}
Given a morphism $\Phi\col A\to B$ of simplicial $R$-algebras, a simplicial resolution of
$B$ over $A$ exists. As noted before, see \ref{resolutions:uniqueness}, such a
resolution is unique up to homotopy of simplicial $A$-algebras.

The following result is clear from the properties of the construction in \ref{killingcycles}.

\begin{Proposition}
  When $R$ is noetherian and the $R$-algebra $S$ is finitely generated, $S$ admits a
  simplicial resolution $\free RX$ with $\card(X_n)$ finite for each $n$. \qed
\end{Proposition}
\end{itchunk}

Next I will describe an explicit resolution of the homomorphism $R[y]\tra R$ where
$y\mapsto 0$.  It serves both as an elementary example of a simplicial resolution, and 
as a way to construct resolutions of complete intersections; see
\ref{resolutions:hypersurfaces}.

\begin{construction}
\label{resolutions:retracts}
Let $R[y]$ be the polynomial ring over $R$, in the variable $y$, and let $\lambda\col
R[y]\to R$ be the homomorphism of $R$-algebras defined by $\lambda(y)=0$.

What is described below is the bar construction of the $R$-algebra $R[y]$ with
coefficients in $R$; see \cite[Chapter X, \S2]{Mc}. It is built as follows: For each
$n\ge0$, set
\[
B_n = R[y]\otimes_R R[y]^{\otimes n} 
\]  
It is convenient to  write $b[b_1|\cdots|b_n]$ for the element $b\otimes b_1\otimes \cdots
\otimes b_n$ in $B_n$.

Evidently, $B_n$ is a polynomial algebra over $R[y]$ over the set of $n$ indeterminates
$X_n=\{x_{nk}\}_{k=0}^{n-1}$, where
\[
x_{nk}= [1|\cdots|y|\cdots|1] \quad \text{with $y$ the $(k+1)$st tensor.}
\]
For each $0\leq i\leq n$, extend the mapping
\[
d_i([b_1|\dots|b_n]) = \begin{cases} 
b_1[b_2|\dots|b_n] & \text{for $i=0$} \\
[b_1|\dots|b_ib_{i+1}|\dots |b_n]& \text{for $1\leq i\leq n-1$}\\
[b_1|\dots|b_{n-1}]\lambda(b_n) & \text{for $i=n$}
\end{cases} 
\]
to a homomorphism of $R[y]$-algebras $d_i\col B_n \to B_{n-1}$. In the same vein,
for each $0\leq j\leq n$ extend the mapping
\[
s_j([b_1|\dots|b_n]) = \begin{cases} 
[1|b_1|\dots|b_n] & \text{for $j=0$} \\
[b_1|1|\dots|b_jb_{j+1}|\dots |b_n]& \text{for $1\leq j\leq n-1$}\\
[b_1|\dots|b_n|1] & \text{for $j=n$}
\end{cases} 
\]
to a homomorphism of $R[y]$-algebras $s_j\col B_n \to B_{n+1}$.

\begin{subexercise}
The $R[y]$-algebra $B=\{B_n\}$, with face and degeneracy operators defined above, is a free
simplicial extension of $R[y]$.
\end{subexercise}

Note that the homomorphism $\lambda$ consists of killing the cycle $y$ in $R[y]$.

\begin{subexercise}
  Reconcile the construction given in this paragraph with the free extension $\free
  {R[y]}{\{x\}\mid \dd(x)=y}$ obtained from \ref{killingcycles}.
\end{subexercise}
\end{construction}

The relevant properties of this free extension are as follows.

\begin{lemma}
\label{bar:koszul}
The canonical surjection $\eps \col B\to R$ is a weak equivalence, and hence
a simplicial resolution of $R$ over $R[y]$.

Let $K$ denote the complex of $R[y]$-module: $0\to R[y]\xra{y}R[y]\to 0$ concentrated in
degrees $0$ and $1$. The morphism of complexes $\nu\col K\to B$ defined by $\nu_0 =
\id^{R[y]}$ and $\nu_1(r) = r\otimes y$, is a homotopy equivalence.
\end{lemma}
\begin{proof}[Sketch of a proof]
  One way to prove this result is as follows: It can be checked directly that $\nu$ is
  compatible with the surjections $\eps\col B\to R$ and $\pi\col K\to R$. Moreover, both
  $\eps$ and $\pi$ are quasi-isomorphisms; this is clear for $\pi$, and is not hard to
  prove for $B$, see \cite[Chapter X, \S2]{Mc}. Thus, since $B$ and $K$ are both bounded
  below complexes of free $R[y]$-modules, it follows that $\nu$ is a homotopy equivalence.

  Another method is to prove first that $\nu$ is a homotopy equivalence, and so deduce
  that $\eps$ is a weak equivalence. I will leave it to you to construct a homotopy
  inverse to $\nu$.
\end{proof}

Given the preceding construction, it is easy to get a simplicial resolution of $R$ over
$R[y_1,\dots,y_d]$. The underlying idea is explained in the following exercise.

\begin{exercise}
\label{resolutions:tensors}
Let $K$ be a commutative ring and let $R'\to S'$ and $R''\to S''$ be homomorphism of
commutative $K$-algebras, such that $R'$ and $S'$ are  flat as $K$-modules.

  If $B'$ and $B''$ are simplicial resolutions of $S'$ over $R'$, and of $S''$ over $R''$,
  respectively, then $B'\otimes_K B''$ is a simplicial resolution of $S'\otimes_KS''$ over
  $R'\otimes_KR''$.
\end{exercise}

Building on the Construction \ref{resolutions:retracts}, I will
describe a simplicial resolution of the $R$-algebra $R/(r)$, when $r$
is a nonzerodivisor on $R$.

\begin{construction}
\label{resolutions:hypersurfaces}
Let $S=R/(r)$, and $\vf\col R\to S$ the canonical surjection.  

Let $R[y]$ denote the polynomial ring over $R$, in the variable $y$, and $\lambda\col
R[y]\to R$ the homomorphism defined by $\lambda(y)=0$. Consider the diagram of
homomorphisms of commutative rings:
\[
\xymatrixrowsep{2.5pc}
\xymatrixcolsep{2.5pc}
\xymatrix{
R[y] \ar@{->}[r]^{\lambda} & R \ar@{->}[d]\\
R  \ar@{<-}[u]^{\psi}    & S=R^{\psi}\otimes_{R[y]}{}^\eps\! R \ar@{<-}[l]_-{\vf}}
\]
Let $R[y]\to B\to R$ be the simplicial resolution of $\lambda$; see Construction
\ref{resolutions:retracts} and Lemma \ref{bar:koszul}. Set
\[
A= R^{\psi}\otimes_{R[y]}B\,.
\]
This is a free extension of $R$.  Base change along along $\psi$ yields a diagram
$R\to A \to S$ of morphisms of simplicial algebras. Since the complex of $R[y]$-modules
$B$ is homotopy equivalent to the complex $K$, defined in Lemma \ref{bar:koszul}, one
obtains that the complex of $R$-modules $A$ is homotopy equivalent to the complex
$R^{\psi}\otimes_{R[y]}K$, that is to say, to the complex:
\[
\xymatrixrowsep{2.5pc}
\xymatrixcolsep{2.5pc}
\xymatrix{
0 \ar@{->}[r] & R \ar@{->}[r]^{r} & R \ar@{->}[r] & 0}
\]
Therefore, we conclude:
\[
\htpy nA=
\begin{cases} 
S & \text{if $n=0$} \\
(0\col_R r) & \text{if $n=1$} \\
0 & \text{if $n\ge 2$} 
\end{cases}
\]
Given this calculation, the proof of the following result is clear.

\begin{Lemma}
  If $r$ is a nonzerodivisor on $R$, then $A$ is a simplicial resolution of the
  $R$-algebra $S$. \qed
\end{Lemma}

The modules and face and degeneracy maps in $A$ are described completely by the
data in \ref{resolutions:retracts}.
\end{construction}

\begin{exercise}
  Suppose that $r_1,\dots,r_c$ is an $R$-regular sequence, see \ref{regular:sequences}.
  Mimic the proof of the preceding lemma to construct a simplicial resolution of $R/(\bsr)$
  over $R$. Hint: use Exercise \ref{resolutions:tensors} and Remark \ref{regular:koszul}.
\end{exercise}

\section{The cotangent  complex}
\label{The cotangent complex}
We are now prepared to introduce the protagonist of these notes: the cotangent complex of
a homomorphism.  This section describes one construction of the cotangent complex, and a
discussion of its basic properties.

\begin{itchunk}{K\"ahler differentials.}
  Since $\kahler -R$ is a functor on the category of $R$-algebras, it extends to a functor
  on the category of simplicial $R$-algebras: Given a simplicial $R$-algebra $A$,
one obtains a simplicial $A$-module $\kahler AR$, with
\[
(\kahler AR)_n = \kahler {A_n}R 
\]
for each $n$, and face maps and degeneracies induced by those on $A$. Moreover, each
morphism $\Phi\col A\to B$ of simplicial $R$-algebras induces a morphism of simplicial
$R$-modules $\kahler \Phi R \col \kahler AR \to \kahler BR$. All this is clear from the
properties of $\kahler -R$ discussed in Section \ref{Kahler differentials}.
\end{itchunk}

\begin{itchunk}{The cotangent complex.}
\label{cotangentcomplex}
Let $\vf\col R\to S$ be a homomorphism of commutative rings.
Let $A$ be a simplicial resolution of $S$ over $R$, and set
\[
\sctan \vf = \kahler AR \otimes_A S\,.
\]
Thus, $\sctan \vf$ is a simplicial $S$-module; the associated complex of $S$-modules is
called the \emph{cotangent complex} of $S$ over $R$; more precisely, of $\vf$.
This too we denote $\sctan\vf$.

Simplicial resolutions are unique up to homotopy, see \ref{resolutionsexist}, and
$\kahler -R$ transforms homotopy equivalent morphisms of simplicial algebras into homotopy
equivalent morphisms of simplicial modules, so the complex $\sctan \vf$ is well defined in
the homotopy category of complexes of $S$-modules; this is explained in the next
paragraph.  It is in this sense that one speaks of \emph{the} cotangent complex.
\end{itchunk}

The crucial point is the following.

\begin{itchunk}{Weak equivalences and differentials.}
\label{kahler:we}
If a morphism of free simplicial $R$-algebras $\Phi\col \free RX \to \free RY$ is a weak
equivalence, then the induced morphism of simplicial $R$-modules $\kahler {\free RX}R\to
\kahler {\free RY}R$ is also a weak equivalence.

The idea is that $\Phi$ admits a homotopy inverse, and hence there is a homotopy inverse
also to $\kahler {\Phi}R$.  Perhaps the best way to formalize this argument is to use the
model category structures on the categories of simplicial $R$-algebras and on simplicial
$R$-modules, see \cite[\S1]{Qu:ams}.
\end{itchunk}

\begin{itchunk}{Homotopies.}
\label{kahler:homotopies}
If $\Phi,\Psi\col \free RX\to B$ are homotopic morphisms of simplicial $R$-algebras, then the
induced morphisms of simplicial $R$-modules $\kahler {\Phi}R$ and $\kahler {\Psi}R$, from
$\kahler {R[X]}R$ to $\kahler BR$, are homotopic.

Indeed, suppose $\free R{X,X,Y}$ is a cylinder object for the $R$-algebra $\free RX$, see
\ref{homotopy}.  Thus, there is diagram of simplicial $R$-algebras
\[
\xymatrixcolsep{2pc}
\xymatrix{
\free R{X,X} \ar@{>->}[r] & \free R{X,X,Y} \ar@{->}[r]^-{\simeq} & \free RX}
\]
where the composed is the product map.  Applying $\kahler -R$ yields a diagram of
simplicial $R$-modules
\[
\xymatrixcolsep{2pc}
\xymatrix{
 \kahler {\free R{X,X}}R \ar@{>->}[r] & \kahler{\free R{X,X,Y}} R \ar@{->}[r]^{\simeq}
  &\kahler{\free RX}R}
\]
The crucial information in the diagram is that the arrow on the right is a weak
equivalence; this is by \ref{kahler:we}.  Concatenating this diagram with the natural
morphism of simplicial $R$-modules
\[
\xymatrixcolsep{1.6pc}
\xymatrix{
\kahler{\free RX}R \bigoplus \kahler{\free RX}R
  \ar@{>->}[r] &
\big(\free RX\otimes_R\kahler{\free RX}R\big)\bigoplus\big(\kahler{\free RX}R\otimes_R\free RX\big) 
 \cong \kahler {\free R{X,X}}R}
\]
where the isomorphism is by Exercise \ref{kahler:products}, one obtains
a diagram
\[
\xymatrixcolsep{2pc}
\xymatrix{
\kahler{\free RX}R \bigoplus \kahler{\free RX}R
  \ar@{>->}[r] & \kahler{\free R{X,X,Y}} R \ar@{->}[r]^{\simeq}
  &\kahler{\free RX}R}
\]
of simplicial $R$-modules. It is easy to check the composed is the canonical map: $
(a,b)\mapsto a + b$.  The diagram above is tantamount to the statement that $\kahler
{\free R{X,X,Y}}R$ is a cylinder object for the simplicial $R$-module $\kahler
{\free RX}R$, see \cite[(2.4)]{Pg:uc}.

Now, applying $\kahler -R$ to the diagram defining a homotopy between $\Phi$ and $\Psi$,
see \ref{homotopy}, one obtains a commutative diagram of simplicial $R$-modules
\[
\xymatrixrowsep{2.5pc}
\xymatrixcolsep{2.5pc}
\xymatrix{
\kahler {\free RX}R \bigoplus \kahler {\free RX}R \ar@{->}[dr]_{\kahler{\Phi}R + \kahler {\Psi}R}
  \ar@{->}[r]  & \kahler{\free R{X,X,Y}}R \ar@{->}[d]^{\simeq} \\
  &\kahler BR}
\]
Since $\kahler{\free R{X,X,Y}}R$ is a cylinder object for $\kahler {\free RX}R$, the
diagram above means that $\kahler{\Phi}R$ and $\kahler{\Psi}R$ are homotopic, see
\cite[(2.7)]{Pg:uc}.
\end{itchunk}

\begin{exercise}
\label{ctan:indecomposables}
Let $A\xra{\eps} S$ be a simplicial resolution, as above. Then $A\otimes_RS$ is a
simplicial $R$-algebra. Let $J$ be the kernel of the morphism of simplicial $S$-algebras
$A\otimes_RS \to \simp{S}$, where $\eps(a\otimes s)=\eps(a)s$. Note that $J$ is a
simplicial ideal in $A\otimes_RS$; that is to say, $J$ is a simplicial
$(A\otimes_RS)$-submodule of $A\otimes_RS$.

  Prove that one has an isomorphism of simplicial $S$-modules:
\[
\sctan \vf \cong J/J^2\,.
\]  
Here $J^2$ is the simplicial ideal in $A\otimes_RS$ with $(J^2)_n= (J_n)^2$.
\end{exercise}

\begin{vista}
  The gist of the preceding exercise is that one may view the cotangent complex as
  `derived indecomposables'.  There are other interpretations of the cotangent complex: as
  the derived functor of the abelianization functor, see \cite[(4.24)]{Pg:uc}; as cotriple
  homology, see \cite[\S8.8]{Wi}.

  In \cite{An:sm}, Andr\'e introduces cotangent complex as in \ref{cotangentcomplex},
  but by using a canonical resolution of the $R$-algebra $S$. This has the benefit that
  one does have to worry that that it is well-defined, and so avoids, in particular, the
  discussion in \ref{kahler:homotopies}. However, to establish any substantial property
  of cotangent complexes, one would have to prove that they can be obtained from any
  simplicial resolution, and so he does.
\end{vista}

\begin{remark}
\label{ctan:freecomplex}
Let $A=\free RX$ be a simplicial resolution of the $R$-algebra $S$.  For each integer $n$,
one has $A_n=R[X_n]$, so Exercise \ref{kahler:freeextension} yields: $\kahler {A_n}R$ is
a free $A_n$-module, and hence $(\sctan \vf)_n$ is a free $S$-module, on a basis of
cardinality $\card(X_n)$.
\end{remark}  

\begin{itchunk}{Andr\'e-Quillen homology and cohomology.}
  The cotangent complex of $\vf$ is well-defined up to homotopy of complexes of
  $S$-modules, so for each $S$-module $N$ and integer $n$, the following $S$-modules are
  well-defined:
\begin{gather*}
\aq nSRN = \HH n{\sctan \vf\otimes_SN} \quad\text{and}\quad
\aqc nSRN = \HH {-n}{\Hom S{\sctan \vf}N} 
\end{gather*}
These are the $n$th \emph{Andr\'e-Quillen homology}, respectively, \emph{Andr\'e-Quillen
  cohomology}, of $S$ over $R$ with coefficients in $N$.

The cotangent complex is a complex of free $S$-modules concentrated in non-negative
degrees, therefore
\[
\aq nSRN = \Tor nS{\sctan \vf}N \quad\text{and}\quad
\aqc nSRN = \Ext nS{\sctan \vf}N\,.
\]
Given this interpretation, a standard argument in the homological algebra of complexes,
see, for instance, \cite[(2.4P), (2.4F)]{AF}, yields the result below.  For any
complex $L$ of $S$-modules, $\fd_SL$ is the \emph{flat dimension} of $L$; thus,
$\fd_SL\leq n$ means that $L$ is quasi-isomorphic to a complex $0\to F_n\to \cdots \to
F_i\to 0$ of flat $S$-modules. The number $\pd_SL$ is the \emph{projective dimension} of
$L$.

\begin{proposition}
\label{aq:fdim}
Let $n$ be a non-negative integer. 

One has $\aq iSR-=0$ for $i\geq n+1$ if and only if $\fd_S(\sctan \vf)\leq n$.

One has $\aqc iSR-=0$ for $i\geq n+1$ if and only if $\pd_S(\sctan \vf)\leq n$. \qed
\end{proposition}
\end{itchunk}

Next I describe the cotangent complex in two cases of interest; it turns out that these
are essentially the only contexts in which one has information in closed form on the
cotangent complex.

\begin{proposition}
\label{ctan:polynomialring}
If $S=R[Y]$, a polynomial ring over $R$ in variables $Y$, then 
the $S$-module $\kahler SR$ is free, and
\[
\sctan\vf \simeq \kahler SR
\]
as complexes of $S$-modules. Thus, $\aq nSRN =0 = \aqc nSRN$ for $n\ge 1$.
\end{proposition}

\begin{proof}
  The freeness of $\kahler SR$ is the content of Exercise \ref{kahler:freeextension}.
  Note that the simplicial $R$-algebra $\simp{S}$ is itself a simplicial resolution of $S$
  over $R$. Therefore, one has the first isomorphism below
\[
\sctan \vf \cong \kahler {\simp{S}}R\otimes_{\simp{S}}S \cong \simp{\kahler SR}\,.
\]
The second isomorphism is verified by inspection.  Thus, as a complex of $S$-modules
$\sctan \vf$, is isomorphic to
\[
\cdots \lra \kahler SR \xra{0} \kahler SR \xra{1}\kahler SR \xra{0} \kahler SR \to 0\,.
\]
Hence, $\sctan\vf$ is homotopy equivalent to $\kahler SR$. The remaining assertions 
now follow, since the $S$-module $\kahler SR$ is free.
\end{proof}

\begin{proposition}
\label{ctan:hypersurface}
If $S=R/(r)$, where $r$ is a nonzerodivisor in $R$, and $\vf\col R\to S$ is the
surjection, then
\[
\sctan \vf \simeq \susp S\,.
\]
\end{proposition}

\begin{proof}
  The proof uses the notation in \ref{resolutions:retracts} and
  \ref{resolutions:hypersurfaces}.

  It is clear that $\sctan \vf$, which equals $\kahler AR\otimes_RS$, is a complex of free
  $S$-modules beginning in degree $1$, and with
\[
(\sctan \vf)_n = \bigoplus_{i=0}^{n-1}Sx_{ni} \qquad \text{for each $n\ge 1$.}
\]
In describing the differential on $\sctan \vf$, it is useful to introduce the following
symbol: for each integer $n$, set
\[
\eps(l,m) = \sum_{k=l}^m (-1)^k = \begin{cases}
 0 & \text{if $m-l$ is even;} \\ 
-1 & \text{if $m-l$ is odd.}
\end{cases}
\]
With this notation, using the description of $A$ ensuing from \ref{resolutions:retracts}
and Exercise \ref{kahler:freeextensionmaps}, one finds that the differential on $(\sctan
\vf)_n$ is given by
\[
\dd(x_{ni}) =  
\begin{cases} 
\eps(1,n) x_{n-1,0} & \text{for $i=0$;}\\
\eps(0,i) x_{n-1,i-1} + \eps(i+1,n)x_{n-1,i} & \text{for $1\leq i\leq n-2$;}\\
\eps(0,n-1) x_{n-1,n-2} & \text{for $i=n-1$.}
\end{cases}
\] 
The entries of the matrix representing the differential 
\[
\dd_n\col (\sctan \vf)_n \to (\sctan \vf)_{n-1}
\]
are either $0$ or $1$.  I claim that $\sctan\vf$ is homotopy equivalent to $\susp S$.

Indeed, since the matrices representing the differentials consist of
zeros and ones, it suffices to verify this assertion when $S=\BZ$
(why?); in particular, we may assume $S$ is noetherian.  A routine
calculation establishes that for any homomorphism $S\to l$, where $l$
is a field, one has
\[
\rank_l(\dd_n\otimes_Sl) = \begin{cases} 
\frac{n-2}2 & \text{if $n\geq 2$ is even;} \\ 
\frac{n+1}2 & \text{if $n\geq 3$ is odd.} 
\end{cases}
\]
Therefore, for each integer $n\geq 2$, one has that
\[
n= \rank_l(\dd_n\otimes_Sl) + \rank_l(\dd_{n+1}\otimes_Sl)
\]
It now remains to do Exercise \ref{exactness:test} below, noting that $\HH 1{\sctan\vf}=S$.
\end{proof}

\begin{exercise}
\label{exactness:test} Let $S$ be a noetherian ring and $L=\{L_n\}_{n\ges 1}$ a complex of finite
  free $S$-modules such that for each prime ideal $\fq$ in $S$, one
  has
\[
\rank_S(L_n)= \rank_l(\dd_n\otimes_Sk(\fq)) + \rank_l(\dd_{n+1}\otimes_Sk(\fq))
\quad\text{for $n\geq 2$.}
\]
Prove that $\HH nL=0$ for $n\geq 2$, and that $L$ is homotopy equivalent to $\susp \HH 1L$.
\end{exercise}

\begin{vista}
  With better machinery one can give more efficient proofs of Proposition
  \ref{ctan:hypersurface}. For instance, writing $\Delta^1$ for the standard 1-simplex, it
  is easy to verify that $\sctan \vf$ is the free the simplicial $S$-module on the
  simplicial set $\Delta^1/\partial \Delta^1$, which implies the desired statement about
  its homotopy, see \cite[(1.15)]{Pg:uc}.  

  Alternatively, one could note that $\sctan\vf$ is the simplicial complex corresponding
  to the chain complex with $S$ in degree 1 (and so zero differential) under the Dold-Kan
  correspondence~\cite[(4.1)]{Pg:uc}, so its homotopy is $S$.  

  Exercise \ref{ctan:polyringrevisited} outlines a third approach. The argument presented
  above was intended to show that, sometimes, one can work directly with simplicial
  resolutions and compute cotangent complexes. Unfortunately, this is perhaps the only
  instance when this is possible, unless one is in characteristic zero, see
  \cite[(9.5)]{Qu:ams}.
\end{vista}

\section{Basic properties}
\label{Basic properties}

This section is a precis of basic properties of the cotangent complex; usually, they are
accompanied by corresponding statements concerning the Andr\'e-Quillen homology modules.
The analogues for cohomology, which are easy to guess, are generally omitted.

As before, let $\vf\col R\to S$ be a homomorphism of rings.

\begin{itchunk}{Functoriality.}
\label{functoriality}
The functor $\sctan \vf\otimes_S- $, defined on the homotopy category of complexes of
$S$-modules, is exact. Therefore, the sequence $\{\aq nSR-\}_{n\in\BZ}$ is a homological
functor on the category of $S$-modules.
\end{itchunk}

\begin{itchunk}{Normalization.}
\label{normalization}
There are isomorphisms of functors
\[
\aq0SR- \cong \Omega_{S\var R}\otimes_S-
\qquad\text{and}\qquad
\aq nSR-=0\quad\text{for each}\quad n<0\,,
\]
where $\Omega_{S\var R}$ denotes the $S$-module of K\"ahler differentials of $S$ over $R$.

Indeed, this is immediate from the right-exactness of $\kahler -R$; see \eqref{jz:kahler}.
\end{itchunk}

\begin{itchunk}{Base change.}
\label{basechange}
Consider a commutative diagram 
\[
\xymatrixrowsep{2.5pc}
\xymatrixcolsep{3.5pc}
\xymatrix{
R'\ar@{->}[r]_{\vf'} \ar@{->}[d]_{\rho} & S' \ar@{->}[d] \\
R\ar@{->}[r]_-{\vf'\otimes_{R'}R=\vf}   & (S'\otimes_{R'}R)\cong S
}
\]
of homomorphisms of rings. It induces a morphism of complexes of $S$-modules:
\[
\sctan {\vf'}\otimes_{R'} R\lra \sctan \vf
\]
which is well defined up to homotopy.  This morphism is a homotopy equivalence when $\Tor
n{R'}{S'}R=0$ for $n\ge 1$; for instance, when either $\vf'$ or $\rho$ is flat.  In this
case, one has isomorphisms of functors
\[ 
\aq n{S}{R}- \cong \aq n{S'}{R'}- \quad\text{for each $n\in \BZ$,}
\]
where $S$-modules are viewed as $S'$-modules via the homomorphism $S'\otimes_{R'}\rho\col
S'\to S$.

Indeed, let $A'\to S'$ be a simplicial resolution of $S'$ over $R'$. This induces a
morphism of simplicial $R$-algebras:
\[
A'\otimes_{R'}R \lra S'\otimes_{R'}{R}=S
\]
Evidently, $A'\otimes_{R'}R$ is a free simplicial extension of $R$.  Thus, if $A$ is a
simplicial resolution of $S$ over $R$, the lifting property yields a morphism of
simplicial $R$-algebra $A'\otimes_{R'}R\to A$, well defined up to homotopy, see
\ref{free:lifting1} and \ref{free:lifting2}.  By functoriality, $\kahler -R$ induces a
morphism of complexes of $S$-modules:
\[
\kahler {(A'\otimes_{R'}R)}{R} \otimes_{(A'\otimes_{R'}R)}S 
   \lra \kahler AR\otimes_A S = \sctan \vf\,,
\]
well-defined up to homotopy of complexes of $S$-modules.  It remains to identify the
complex on the left, and this is accomplished below:
\begin{align*}
\kahler {(A'\otimes_{R'}R)}R \otimes_{(A'\otimes_{R'}R)} S
           &\cong\big(\kahler {A'}{R'}  \otimes_{A'} R\big)
                       \otimes_{(A'\otimes_{R'}R)} (S'\otimes_{R'}R) \\
           &\cong  \big(\kahler{A'}{R'}\otimes_{A'} S'\big) \otimes_{R'}R\\
           &= \sctan {\vf'}\otimes_{R'}R
\end{align*}
The isomorphisms are all verified directly.

Suppose $\Tor n{R'}{S'}R=0$ for $n\ge1$. As noted in \ref{resolution:underlying}, one has
an isomorphism
\[
\htpy n{A'\otimes_{R'}R} \cong \Tor n{R'}SR \quad \text{for each $n$.}
\]
Therefore, the augmentation $A'\otimes_{R'}R\to S$ is a weak equivalence, and hence
$A'\otimes_{R'}R$ is a simplicial resolution of $S$ over $R$. Thus, the morphism
$A'\otimes_{R'}R\to A$ is a homotopy equivalence, and hence so is the induced morphism
$\sctan{\vf'}\otimes_{R'}R\to \sctan \vf$.
\end{itchunk}

Here is a beautiful application, due to Andr\'e, of the preceding property:

\begin{proposition}
\label{ctan:localizationmap}
  Let $U$ be a multiplicatively closed subset of $R$, let $S=U^{-1}R$, and let $\vf\col
  R\to S$ be the localization map. Then $\sctan \vf \simeq 0$.
\end{proposition}

\begin{proof}
  The complex $\sctan \vf$ consists of $S$-modules, and the functor $-\otimes_RS$ is the
  identity on the category of $S$-modules, so one obtains the isomorphism below
\[
\sctan \vf \cong \sctan \vf \otimes_R S \simeq 
  \sctan {\vf\otimes_RS}\simeq \sctan{\id^S}\simeq 0
\]
The first homotopy equivalence holds by Property \ref{basechange}, since the homomorphism
$\vf$ is flat, the second one holds because the homomorphism $\vf\otimes_RS\col
S\otimes_RS\to S$ is the identity, while the last one follows, for example, from
Proposition \ref{ctan:polynomialring}.
\end{proof}

\begin{itchunk}{Localization.}
\label{localization}
Fix a prime ideal $\fq$ in $S$, set $\fp=R\cap \fq$, and denote $\vf_\fq\col R_\fp\to
S_\fq$ the localization of $\vf$ at $\fq$. One has a homotopy equivalence
\[
\sctan {\vf_\fq} \simeq S_\fq\otimes_S \sctan \vf 
\]
of complexes of $S_\fq$-modules.  In particular, for each $n\in \BZ$, there is an
isomorphism of functors of $S$-modules
\[
\aq nSR-_\fq \cong \aq n{S_\fq}R{-_\fq}\cong\aq n{S_\fq}{R_{\fp}}{-_\fq}
\]
See \cite{An:hca} for a proof of these assertions. Alternatively:

\begin{subexercise}
Prove the assertions above.
\end{subexercise}
\end{itchunk}

In this context, one has the following useful remark which permits one to reduce many
problems concerning the vanishing of Andr\'e-Quillen homology to the case of homomorphisms
of local rings.

\begin{proposition}
\label{vanishing:andre}
Let $\vf\col R\to S$ be a homomorphism of rings.  For each integer $n$, the following
conditions are equivalent:
\begin{enumerate}[\quad\rm(a)]
\item 
$\aq nSR-=0$ on the category of $S$-modules;
\item $\aq n{S_\fq}{R_{\fq\cap R}}-=0$ on the category of $S_\fq$-modules, for each
  $\fq\in \Spec S$. \qed
\end{enumerate}
\end{proposition}

The proof of this result is elementary, given property \ref{localization}.  Under an
additional hypothesis on $\vf$, there is a significant improvement to the preceding
result; see Proposition \ref{vanishing:local}.

\begin{itchunk}{Transitivity.} 
\label{transitivity}
Each homomorphism of rings $\psi\col Q\to R$ induces the following exact triangle 
in the homotopy category of complexes of $S$-modules:
\[
(S\otimes_R\sctan \psi) \lra \sctan{\vf\circ\psi}\lra \sctan \vf \lra \susp (S\otimes_R\sctan \psi)
\]
This induces an exact sequence of functors of $S$-modules
\[
\cdots\lra \aq {n+1}SR-\lra\aq nRQ- \lra \aq nSQ- \lra \aq nSR- \lra \cdots\,.   
\]
It is called the \emph{Jacobi-Zariski sequence} arising from the diagram $Q\to R\to S$.
It extends \eqref{jz:kahler} to a long exact sequence of $S$-modules.

For a proof of this assertion, see \cite[(4.32)]{Pg:uc}, or \cite[(5.1)]{Qu:ams}.
\end{itchunk}

Use the transitivity sequence to solve the following exercises.

\begin{exercise}
\label{ctan:localization}
Let $U\subset S$ be multiplicatively closed subset, and  $\eta\col S\to U^{-1}S$ the
localization map. Prove that one has a homotopy equivalence
\[
\sctan {\eta\circ\vf}\simeq U^{-1}S\otimes_S\sctan \vf
\]
of complexes of $U^{-1}S$ modules.
\end{exercise}

\begin{exercise}
\label{retract}
Let $\psi\col S\to R$ be a homomorphism of rings such that the map $\vf\psi\col S\to S$
equals $\id^S$; said otherwise, $S$ is an algebra retract of $R$.

Prove that one has a homotopy equivalence of complexes of $S$-modules:
\[
\sctan \vf \simeq \susp (\sctan \psi \otimes_RS)\,.
\]
In particular, $\aq nSR- \cong \aq {n-1}RS-$ as functors of $S$-modules.
\end{exercise}

\begin{exercise}
\label{ctan:polyringrevisited}
Use the preceding exercise, and the discussion in Construction
\ref{resolutions:hypersurfaces}, to prove Proposition \ref{ctan:hypersurface}.
\end{exercise}

\begin{itchunk}{Finiteness.}
\label{finiteness}
Suppose $R$ is noetherian and $\vf$ is essentially of finite type.

The complex $\sctan \vf$ is then homotopic to a complex
\[
\cdots \lra L_n \lra L_{n-1} \lra \cdots \lra L_1 \lra L_0 \lra 0\,,
\]
where for each $n$, the $S$-module $L_n$ is finitely generated and free.  Thus,
when the $S$-module $N$ is finitely generated so are $\aq nSRN$ and $\aqc nSRN$.

Indeed, by hypothesis $\vf$ admits a factorization
\[
R\xra{\eta} U^{-1}R[Y]\xra{\vf'} S
\]
with $\card(Y)$ finite, $U$ a multiplicatively closed subset of $R[Y]$, and $\vf'$ a
surjective homomorphism of rings. Since $\eta$ factors as $R\to R[U]\to U^{-1}R[Y]$, it
follows from Proposition \ref{ctan:polynomialring} and Exercise \ref{ctan:localization}
that $\sctan{\eta}$ is equivalent to a complex of finitely generated free $U^{-1}R[Y]$
modules.  On the other hand, Proposition \ref{resolutionsexist} and Remark
\ref{ctan:freecomplex} imply that the complex $\sctan{\vf'}$ consists of finitely
generated free modules $S$-modules. Now the desired result is a consequence of
\ref{transitivity}, applied to the diagram above.
\end{itchunk}

\begin{itchunk}{Low degrees.}
\label{low degrees}
As usual, low degree cohomology modules admit alternative interpretations.  First, a piece
of notation: $S\ltimes N$ denotes the ring with being $S\oplus N$ the underlying abelian
group and product given by $(s,x)(t,y)= (st,sy+tx)$.

To begin with, the $S$-module $\aqc 0SRN$, which is $\der RSN$, is the set of $R$-algebra
homomorphisms $\alpha \col S\ltimes N \to S\ltimes N$ extending the identity map both on
$N$ and on $S$, that is to say, such that the following diagram commutes:
\[
\xymatrixcolsep{2pc}
\xymatrixrowsep{2pc}
\xymatrix{
0 \ar@{->}[r] & N \ar@{->}[r] \ar@{=}[d] & S\ltimes N \ar@{->}[r]^{\eps}\ar@{->}[d]^{\alpha} &
            S \ar@{->}[r] \ar@{=}[d] & 0 \\
0 \ar@{->}[r] & N \ar@{->}[r] & S\ltimes N \ar@{->}[r]^{\eps} & S \ar@{->}[r] & 0\,.
}
\]
Here $\eps$ is the canonical surjection.  This claim is not hard to verify; see
\cite[\S25]{Ma}.

The $S$-module $\aqc 1SRN$ is the set of isomorphism classes of extensions
of $R$-modules
\[
\xymatrixcolsep{2pc}
\xymatrixrowsep{2pc}
\xymatrix{
0 \ar@{->}[r] & N \ar@{->}[r]^{\iota} & \wt S \ar@{->}[r]^{\eps} & S \ar@{->}[r] & 0\,,
}
\]
where $\eps$ is a homomorphism of $R$-algebras with $\Ker(\eps)^2=(0)$, and the given $S$
module structure on $N$ coincides with the one induced by $\iota$; see \cite[Chapter
XVI]{An:hca} for a proof.

When $S=R/I$, one has $\aq 0SRN =0 =\aqc 0SRN$, see \ref{normalization},
and
\begin{equation}
\label{surjective:degree1}
\aq 1SRN = (I/I^2) \otimes_S N \quad\text{and}\quad \aqc 1SRN =\Hom S{I/I^2}N
\end{equation}
These claims are justified by Proposition \ref{aqvstor:degree1}.
\end{itchunk}

\section{Andr\'e-Quillen homology and the Tor functor}

In this section we discuss the relationship between the Andr\'e-Quillen homology modules
$\{\aq nSRN\}$, where $N$ is an $S$-module, and the $S$-modules $\{\Tor nRSN\}$.

Let $\eps\col A\to S$ be a simplicial resolution of the $R$-algebra $S$, and let $J$
denote the simplicial ideal $\Ker(A\otimes_RS\to \simp S)$; see Exercise
\ref{ctan:indecomposables}. 

One has an exact sequence of simplicial $S$-modules
\[
0\to J\to A\otimes_R S\to \simp{S}\to 0\,.
\]
Since $\simp{S}_n=S$ for each $n$, applying $-\otimes_SN$ preserves the exactness of the
sequence above, so passing to the homology long exact sequence yields 
\[
\htpy n{J\otimes_SN} = \begin{cases}
\Ker(S\otimes_RN\to N)& \text{when $n=0$;} \\
\Tor nRSN & \text{when $n\ge 1$.}
\end{cases}
\]
Here one is using Remark \ref{resolution:underlying}.  

The morphism $J\to J/J^2$ induces a morphism $J\otimes_SN\to (J/J^2)\otimes_SN$ of
simplicial modules. In homology this yields, keeping in mind Exercise
\ref{ctan:indecomposables} and the preceding display,  homomorphisms
of $S$-modules:
\[
  \Tor nRSN \lra \aq nSRN \quad \text{for $n\ge 1$.}
\]
Naturally, the properties of this map are determined by those of the simplicial ideal $J$,
which in turn reflects properties of the $R$-algebra structure of $S$. The following
result, which justifies the claim in \eqref{surjective:degree1}, is one manifestation of
this phenomenon.

\begin{proposition}
\label{aqvstor:degree1}
Assume that $\vf$ is surjective, and set $I=\Ker(\vf)$. One has natural isomorphisms of
$S$-modules
\[
\aq 1SRN\cong \Tor 1RSN\cong (I/I^2)\otimes_SN\,.
\]
\end{proposition}

\begin{proof}
  Since $\vf$ is surjective, one may choose a simplicial resolution $A$ of $S$ with
  $A_0=R$. Set $B=A\otimes_RS$. The crucial point in the proof is the following

\emph{Claim}. $\HH 0{J^2}=0=\HH 1{J^2}$.

Indeed, by choice of $A$, one has $J_0=0$, which explains the first equality.  Moreover,
the cycles in $\norm {J^2}_1$ equal $J_1^2$, and hence a sum of elements of the form $xy$,
where $x$ and $y$ are in $J_1$. However, $xy$ is a boundary: the element 
\[
w = s_0(xy) - s_0(x)s_1(y)
\]
is an element in $\norm{J^2}$ and $d_0(w) = xy$. Thus, $\HH 1{J^2}=0$.

Now, in the exact sequence $0\to J^2\to J\to J/J^2\to 0$ of simplicial modules, for each
integer $n$, the $S$-module $(J/J^2)_n$ is free, so one obtains an exact sequence
\[
0\to J^2\otimes_SN \to J\otimes_SN \to (J/J^2)\otimes_S N\to 0\,.
\]
Passing to homology and applying the claim above yields the first of the desired isomorphisms.

As to the second one: $\Tor 1RSN\cong (I/I^2)\otimes_SN$, consider the exact sequence
\[
0 \lra I \lra R \lra S\lra 0
\]
and apply to it the functor $-\otimes_RN$.
\end{proof}

The next theorem was proved by Quillen, see \cite[(6.12)]{Qu:ams}, \cite[Chapter XX,
(24)]{An:hca}; it extends the proposition above, for when $\vf$ is surjective, the
multiplication map $\mu^S_R\col S\otimes_RS\to S$ is bijective.

\begin{theorem}
\label{connectivity}
If $\mu^S_R$ is bijective, then $\HH i{J^n}=0$ for $n\ge 1$ and $i\leq n-1$.\qed
\end{theorem}

This result is a critical component in proving the convergence of a
spectral sequence relating Andr\'e-Quillen homology and the Tor functor:

\begin{itchunk}{The fundamental spectral sequence.}
  Suppose that $\mu^S_R$ is bijective. The $S$-modules underlying the sub-quotients of the
  filtration $\cdots \subseteq J^2\subseteq J\subseteq A$ are free, so one
  obtains a filtration
\[
\cdots \subseteq (J^2\otimes_SN) \subseteq (J\otimes_SN)\subseteq (A\otimes_SN)\,.
\]
This induces a spectral sequence with
\[
\EH 1pq = \left((J^q/J^{q+1})\otimes_SN\right)_{p+q}
\]
and abutting to $\HH{p+q}{A\otimes_SN}=\Tor{p+q}RSN$.  It follows from the connectivity
theorem \ref{connectivity} that
\[
\EH 2pq = \htpy{p+q}{ \big(J^q/J^{q+1}\big)\otimes_SN } =0 \quad \text{for $p\leq -1$.}
\]
Thus, the spectral sequence converges. Given Exercise \ref{ctan:indecomposables}, the
5-term exact sequence arising from the edge homomorphisms of the spectral sequence yields

\begin{proposition}
\label{fss:5term}
If $\mu_R^S$ is surjective, there is an exact sequence of $S$-modules
\begin{align*}
\Tor 3RSN\to \aq3SRN \to &\wedge^2_S\Tor 1RSS\otimes_SN \to \cdots \\
\cdots  &\to \Tor 2RSN\to \aq2SRN\to 0\,.
\end{align*}
\end{proposition}
This result will be used in the study of homomorphisms of noetherian rings, which is the
topic of the next section.
\end{itchunk}

\section{Locally complete intersection homomorphisms}
The remainder of this article concerns the role of Andr\'e-Quillen homology in the study
of homomorphisms of commutative rings. The section focuses on complete intersection
homomorphisms, while the next one is dedicated to regular homomorphisms.  Henceforth, the
tacit assumption is that rings are noetherian. Recently, I learned of a new book by
Majadas and Rodicio \cite{MR} aimed at providing a comprehensive treatment of the basic
results in this topic.

\begin{itchunk}{Regular sequences.}
  \label{regular:sequences} A sequence $\bsr=r_1,\dots,r_c$ of elements of $R$ is said to
  be \emph{regular} if $(\bsr)\ne R$ and $r_i$ is a nonzerodivisor on
  $R/(r_1,\dots,r_{i-1})$ for $i=1,\dots,c$.
\end{itchunk}

For example, in the ring $R[y_1,\dots,y_c]$, the sequence
$y_1,\dots,y_c$ is regular.

\begin{remark}
\label{regular:koszul}
Given an element $r$ in $R$, write $\koszul rR$ for the complex of $R$-modules $0\to
R\xra{r} R\to 0$, with non-zero modules situated in degrees $0$ and $1$. Given a sequence
of elements $\bsr=r_1,\dots,r_c$ in $R$, set
\[
\koszul \bsr R=\, \koszul {r_1}R\otimes_R \cdots \otimes_R \koszul {r_c}R
\]
This is the \emph{Koszul complex} on the elements $\bsr$.

Koszul complexes were applied to the study of regular sequences by Auslander and Buchsbaum
who proved: if $\bsr$ is a regular sequence, then $\HH n{\koszul \bsr R}=0$ for $n\ge 1$,
so the augmentation $\koszul \bsr R\to R/\bsr R$ is a quasi-isomorphism. The converse
holds when $\bsr$ is contained in the Jacobson radical of $R$, see \cite[(16.5)]{Ma}.
\end{remark}

\begin{itchunk}{Locally complete intersection homomorphisms.}
\label{lci:definition}  
Let $\vf\col R\to S$ be a homomorphism of noetherian rings.

When $\vf$ is surjective, it is \emph{complete intersection} if the ideal $\Ker(\vf)$ is
generated by a regular sequence; it is \emph{locally complete intersection} if for each
prime ideal $\fq$ in $S$, the homomorphism $\vf_\fq\col R_{\fq\cap R}\to S_\fq$ is
complete intersection.

When $\vf$ is a homomorphism essentially of finite type, it is \emph{locally complete
  intersection} if in some factorization
\[
R\to R'\xra{\vf'} S
\]
of $\vf$ where $R'$ is of the form $U^{-1}R[Y]$, where $U$ is a multiplicatively closed
subset in $R[Y]$, and $\vf'$ is surjective, the homomorphism $\vf'$ is locally complete
intersection. It is not too hard that this property is independent of the chosen
factorization; it becomes easy, once Theorem \ref{ci:aqcriteria} is on hand.

Avramov has introduced a notion of a complete intersection homomorphism \emph{at a prime
  $\fq$ in $\Spec S$}, and of locally complete intersection homomorphisms, for arbitrary
homomorphisms of noetherian rings. It is based on the theory of `Cohen factorizations';
see \cite[\S1]{Av:am}.

Vanishing of Andr\'e-Quillen homology is linked to the locally complete intersection
property by following result, which was proved by Lichtenbaum and Schlessinger, Andr\'e,
and Quillen in the case when $\vf$ is essentially of finite type, and by Avramov in the
general case.
\end{itchunk}

\begin{theorem}
\label{ci:aqcriteria}
Let $\vf\col R\to S$ be a homomorphism of noetherian rings.

The following conditions are equivalent.
\begin{enumerate}[{\quad\rm(a)}]
\item
$\vf\col R\to S$ is locally complete intersection.
\item
$\aq nSR-=0$ for $n\geq 2$\,.
\item
$\aq 2SR-=0$\,.
\end{enumerate}
\end{theorem}
Condition (b) may be restated as: $\fd_S\sctan \vf\leq 1$; see Proposition \ref{aq:fdim}.

We prove the theorem above when $\vf$ is essentially of finite type. In that case, the
implication (c) $\implies$ (a) is reduced to the more general result below.

\begin{theorem}
\label{ci:localcase}
Let $\vf\col (R,\fm,k)\to (S,\fn,l)$ be a local homomorphism, essentially of finite type.
If $\aq 2SRl=0$, then  $\vf$ is locally complete intersection.
\end{theorem}

\begin{remark}
  \label{cimap:completion} The hypothesis that $\vf$ is essentially of finite type is
  needed: the local homomorphism $\zeta\col (R,\fm,k) \to (\wh R, \fm \wh R,k)$, where
  $\wh R$ is the $\fm$-adic completion of $R$, is flat, so base change along $R\to k$
  yields, by \ref{basechange}, the isomorphism below:
\[
\aq  n{\wh R}Rk \cong \aq nkkk =0 \quad \text{for each $n$.}
\]
However, $\zeta$ is locally complete intersection if and only if the formal fibers of $R$,
that is to say, the fibres of the homomorphism $R\to \wh R$, are locally complete
intersection rings, in the sense of \ref{cirings}, and this is not always the case; see
\cite{Mar} and \cite{Tab}.
\end{remark}

\begin{proof}[Proof of Theorem \ref{ci:localcase}]
  By hypothesis, $\vf$ can be factored as $R\xra{\eta} R'\xra{\vf'} S$, where
  $R'=R[Y]_\fq$, with $\card(Y)$ finite, $\fq$ is a prime ideal in $R[Y]$, and $\vf'$ is a
  surjective homomorphism. Proposition \ref{ctan:polynomialring} and Exercise
  \ref{ctan:localization} yield $\aq n{R'}R-=0$ for $n\ge 2$, so the Jacobi-Zariski
  sequence \ref{transitivity} yields isomorphisms
\[
\aq 2{S}{R'}l \cong \aq 2{S}{R}l=0\,.
\]
Therefore, replacing $R'$ by $R$, one may assume $\vf$ is surjective.
In particular, $k=l$.

Suppose $\Ker(\vf)$ is minimally generated by $\bsr = r_1,\dots,r_c$, so $S=R/(\bsr)R$.
We prove, by an induction on $c$, that the sequence $\bsr$ is regular.  

When $c=1$, so that $S=R/rR$, Proposition \ref{fss:5term}, specialized to $N=k$,
yields an exact sequence
\[
\to \wedge^2 \Tor 1RSS\otimes_Sk \to \Tor 2RSk \to \aq 2SRk \to 0\,.
\]
Note that $\Tor 1RSS = (r)/(r^2)$, so $\Tor 1RSS \otimes_S k \cong k$, and hence
\[
\wedge^2\Tor 1RSS\otimes_Sk \cong \wedge^2 k  = 0\,.
\]
Thus, since $\aq 2SRk=0$, the exact sequence above implies $\Tor 2RSk=0$. The ring $R$ is
local and $R$-module $S$ is finitely generated, so the last equality implies
$\pd_RS\leq 1$, see \cite[\S19, Lemma 1]{Ma}. Since the complex
\[
0\to R\xra r R\to 0
\]
is the beginning of a minimal resolution of $S$, one deduces that it is the minimal
resolution.  In particular, $r$ is a nonzerodivisor on $R$, as required.

Suppose the result has been proved whenever $\Ker(\vf)$ is minimally generated by $c-1$
elements.  Set $R'=R/(r_1,\dots,r_{c-1})R$. The Jacobi-Zariski sequence \ref{transitivity}
arising from the diagram $R\to R'\to S$ yields an exact sequence
\[
\to \aq 2SRk \to \aq 2S{R'}k \to \aq 1{R'}Rk \to \aq 1SRk \to \aq 1S{R'}k \to 0
\]
It follows from Proposition \ref{aqvstor:degree1} that 
\[
\aq 1{R'}Rk\cong k^{c-1}, \quad \aq 1SRk\cong k^c, \quad \text{and }
\aq 1S{R'}k\cong k\,.
\]
Thus, since $\aq 2SRk=0$, the exact sequence above yields an exact
sequence
\[
0 \to \aq 2S{R'}k \to k^{c-1} \to k^c \to k \to 0
\]
Therefore, $\aq 2S{R'}k=0$, and since $S=R'/r_cR'$ the basis of the induction implies
$r_c$ is a nonzerodivisor on $R'$. In particular, $\aq 3S{R'}k=0$, by Proposition
\ref{ctan:hypersurface}.  Given that $\aq 2SRk=0$, the Jacobi-Zariski sequence
\ref{transitivity} yields
\[
\aq 2{R'}Rl \cong \aq 2SRl = 0
\]
Consequently, $\aq 2{R'}Rk=0$. Thus, the induction hypothesis implies the sequence
$r_1,\dots,r_{c-1}$ is regular on $R$. This is as desired, since $r_c$ is regular on
$R'$.
\end{proof}

Here is another simplification which results in the theory of Andr\'e-Quillen homology
when the homomorphism under consideration is essentially of finite type.

\begin{lemma}
\label{vanishing:local}
Let $\vf\col R\to (S,\fn,l)$ be a local homomorphism essentially of finite type.  The
complex of $S$-modules $\sctan \vf$ is homotopic to a complex 
\[
\cdots \lra L_n \lra L_{n-1} \lra \cdots \lra L_1 \lra L_0 \lra 0
\]
of finite free $S$-modules, and with $\dd(L)\subseteq \fn L$.

In particular, for each integer $n$, one has $\rank_S(L_n)=\rank_l\aq nSRl$, so that if
$\aq nSRl=0$, then $\aq nSR-=0$ on the category of $S$-modules.
\end{lemma}

\begin{proof}
  One way to prove this result is to note that, since $\vf$ is essentially of finite type,
  $\sctan \vf$ is homotopy equivalent to a complex $L=\cdots \to L_1\to L_0\to 0$ of
  finite free $S$-modules, see \ref{finiteness}. Since $S$ is local, $L$ is homotopic to
  one with $\dd(L)\subseteq \fn L$; for instance, see \cite[(1.1.2)]{Av:sl}. The desired
  claim is now clear.
\end{proof}

\begin{proof}[Proof of Theorem \emph{\ref{ci:aqcriteria}}]
  We give the argument when $\vf$ is essentially of finite type; see \cite{Av:am} for
  the general case. All three conditions are local properties: condition (a) by
  inspection, and conditions (b) and (c) by Lemma \ref{vanishing:local}. Thus, one
  may assume $\vf\col (R,\fm,k)\to (S,\fn,l)$ is a local homomorphism.

  Now (c) $\implies$ (a) follows from Theorem \ref{ci:localcase}, while (b) $\implies$ (c)
  is obvious.

  (a) $\implies$ (b). Arguing as in the proof of Theorem \ref{ci:localcase}, one may reduce
  to the case where $\vf$ is surjective.  Suppose $\Ker(\vf)$ is minimally generated by
  elements $\bsr = r_1,\dots,r_c$; thus $S=R/\bsr R$.  An elementary induction on $c$,
  using Proposition \ref{ctan:hypersurface} and Property \ref{transitivity}, yields
  $\sctan\vf\simeq \susp S^c$.  Therefore, $\aq nSR-=0$ for $n\geq 2$.
\end{proof}

Now the following exercise should not be too taxing.

\begin{exercise}
  Suppose $\vf$ is essentially of finite type. Prove that when
  $\vf$ is locally complete intersection, in any factorization $R\to
  U^{-1}R[Y]\xra{\vf'} S$ of $\vf$, where $\vf'$ is surjective, the
  homomorphism $\vf'$ is locally complete intersection.
\end{exercise}

Here is an exercise which illustrates the flexibility afforded by the characterization in
Theorem \ref{ci:aqcriteria}. To better appreciate it, try to solve it without using the
machinery of Andr\'e-Quillen homology.

\begin{exercise}
\label{ci:basechange}
Let $\vf\col R\to S$ be a homomorphism of noetherian rings, essentially of finite type,
and let $R\to R'$ be a flat homomorphism.

Prove that if $\vf$ is locally complete intersection, then so is the induced homomorphism
$\vf\otimes_RR'\col R'\to (S\otimes_R R')$, and that the converse holds when $R\to R'$ is
faithfully flat. Hint: for the converse, use the going-down theorem, see
\cite[(9.5)]{Ma}.
\end{exercise}

\begin{itchunk}{Extensions of fields.}
\label{cartier}
Let $\phi\col k\to l$ be a homomorphism of fields.  

\begin{subexercise}
  Prove that when the field $l$ is finitely generated over $k$, the homomorphism $\phi$ is
  locally complete intersection.
\end{subexercise}

It is easy to check that $\phi$ is locally complete intersection in general, in the sense
of \cite{Av:am}.  Thus, $\aq nlk-=0$ for $n\geq 2$, by Theorem \ref{ci:aqcriteria}. The
$l$-vectorspace $\aq 1lkl$ is called the \emph{module of imperfection}, and denoted
$\Gamma_{\!l\mid k}$; see \cite[\S26]{Ma}.

When $h\to k$ is another homomorphism of fields, the Jacobi-Zariski sequence
\ref{transitivity} arising from the diagram $h\to k\to l$ yields an exact sequence of
$l$-vectorspaces:
\[
0\to \Gamma_{\!k\mid h}\otimes_kl \to \Gamma_{l\mid h} \to \Gamma_{l\mid k} \to
\kahler{k}{h}\otimes_kl \to \kahler{l}{h} \to \kahler{l}{k} \to 0
\]
Computing ranks one obtains the Cartier equality, see \cite[(26.10)]{Ma}.
\end{itchunk}

\begin{itchunk}{Regular local rings.}
\label{regularrings}
A local ring $(R,\fm,k)$ is \emph{regular} if the ideal $\fm$ has a set of
generators that form an $R$-regular sequence. This condition translates to: the surjection
$R\to k$ is complete intersection, in the sense of \ref{lci:definition}.  The following
result is a corollary of Theorems \ref{ci:aqcriteria} and \ref{ci:localcase}.

\begin{proposition}
  Let $R$ be a local ring, with residue field $k$. The following conditions are
  equivalent.
\begin{enumerate}[\quad\rm(a)]
\item 
$R$ is regular;
\item
$\aq nkR-=0$ for $n\ge 2$;
\item
$\aq 2kRk=0$. 
\end{enumerate}
When $R$ is regular, $\fm$ its maximal ideal, and $\eps\col R\to k$ is the canonical
surjection, then $\sctan \eps \simeq (\fm/\fm^2)$. \qed
\end{proposition}

This result is a homological characterization of the regularity property, akin to the one
by Auslander, Buchsbaum, and Serre: $R$ is regular iff every $R$-module has finite
projective dimension iff $k$ has finite projective dimension, see \cite[\S19]{Ma}.
\end{itchunk}

\begin{itchunk}{Complete intersections.}
\label{cirings}
Let $(R,\fm,k)$ be a local ring, and let $\wh R$ denote the $\fm$-adic completion of $R$.
Cohen's structure theorem provides a surjection $\eps\col Q\tra \wh R$ with $Q$ a regular
local ring; see \cite[(29.4)]{Ma}. Such a homomorphism $\eps$ is said to be a \emph{Cohen
  presentation} of $\wh R$.

The local ring $R$ is \emph{complete intersection} if in a Cohen presentation $\eps\col
Q\tra \wh R$, the ideal $\Ker(\vf)$ is generated by a regular sequence; that is to say,
$\eps$ is a complete intersection homomorphism.  It is known, and is implicit in the proof
of the result below, that when $R$ is complete intersection, any Cohen presentation of
$\wh R$ is a complete intersection homomorphism.

\begin{proposition}
  Let $R$ be a local ring, with residue field $k$. The following conditions are
  equivalent.
\begin{enumerate}[\quad\rm(a)]
\item
 
$R$ is complete intersection;
\item
$\aq nkR-=0$ for $n\ge 3$;
\item
$\aq 3kRk=0$.
\end{enumerate}
\end{proposition}
\begin{proof}
  Since $k$ is a field, and $\aq nkR-$ commutes with arbitrary direct sums (check this),
  condition (b) is equivalent to:
\begin{enumerate}
\item[\quad\rm(b')]
  $\aq nkRk=0$ for $n\ge 3$.
\end{enumerate}
The homomorphism $R\to \wh R$ is flat, see \cite[(8.8)]{Ma}, so base change
along it yields isomorphisms
\[
\aq nk{\wh R}k \cong \aq nkRk \quad \text{for $n\in\BZ$.}
\]
Therefore, we may assume that $R$ is complete, and hence that there is a surjection
$\eps\col Q\to R$, where $Q$ is a regular local ring, see \ref{cirings}.

Corollary \ref{regularrings} yields $\aq nkQk=0$ for $n\ge 2$, so the Jacobi-Zariski
sequence \ref{transitivity} applied to the diagram $Q\to R\to k$ provides isomorphisms
\[
\aq nkRk \cong \aq{n-1}RQk \quad \text{for $n\ge 2$.}
\]

Now, when $R$ is complete intersection, there is a choice of $\eps$ which is complete
intersection. Then Theorem \ref{ci:aqcriteria} implies $\aq nRQ-=0$ for $n\geq 2$; note
that, since $\eps$ is surjective, we are applying the theorem in a case where it was
proved. Thus, the isomorphisms above imply condition (b').

Conversely, given (c), one obtains $\aq nRQk=0$ for $n=2$, by the displayed isomorphisms.
Now Theorem \ref{ci:localcase} yields that $\eps$ is complete intersection. Hence,
$R$ is complete intersection.
\end{proof}
\end{itchunk}

\begin{exercise}
  Let $(R,\fm,k)$ be a local ring, and $\eps\col Q\tra R$ a surjective homomorphism with
  $Q$ a regular local ring. Prove that the ring $R$ is complete intersection if and only
  if the homomorphism $\eps$ is a complete intersection.
\end{exercise}

\begin{itchunk}{The Quillen conjectures.}
  For homomorphisms of noetherian rings, and essentially of finite type, in
  \cite[(5.6), (5.7)]{Qu:ams} Quillen made the following conjectures:

\medskip

\textbf{Conjecture I.} If $\fd_S\sctan \vf$ and $\fd_RS$ are both finite, then the
homomorphism $\vf$ is locally complete intersection.

\medskip

\textbf{Conjecture II.} If $\fd_S\sctan \vf$ is finite, then $\fd_S\sctan \vf\leq 2$.

\medskip

Recall that $\fd_S\sctan \vf\leq n$ if and only if $\aq iSR-=0$ for $i\geq n+1$, see
Proposition \ref{aq:fdim}, so the Quillen conjectures can be phrased in terms of vanishing of
Andr\'e-Quillen homology functors.

Avramov \cite{Av:am} proved the following result, settling Conjecture I in the
affirmative:

\begin{theorem}
Let $\vf\col R\to S$ be a homomorphism of noetherian rings.

The following conditions are equivalent.
\begin{enumerate}[{\quad\rm(i)}]
\item
$\vf$ is locally complete intersection.
\item
$\aq nSR-=0$ for $n\gg 0$ and $\fd_RS$ is locally finite.
\end{enumerate}
If $S$ has characteristic $0$, then they are also equivalent to
\begin{enumerate}[{\quad\rm(i)}]
\item[{\rm(iii)}]
$\aq mSR-=0$ for some integer $m\ge 2$ and $\fd_RS$ is locally finite. \qed
\end{enumerate}
\end{theorem}

Jim Turner \cite{Jt} gave a different proof of Quillen's conjecture I, in the special case
when $\vf$ is essentially of finite type, and the residue fields of $R$ are all of
positive characteristic.

In \cite{AI:ens}, Conjecture II is settled for homomorphisms that admit algebra retracts:

\begin{theorem}
\label{third}
Let $\vf \col R\to S$ be a homomorphism of noetherian rings such that there exists a
homomorphism $\psi\col S\to R$ with $\vf\col\psi=\id^S$.

The following conditions are equivalent.
\begin{enumerate}[{\quad\rm(i)}]
\item
$\psi_\fp$ is complete intersection for each $\fp\in\Spec R$ with
$\fp\supseteq \Ker(\vf)$\,.
\item
$\aq nSR-=0$ for $n\gg 0$.
\item
$\aq nSR-=0$ for $n\ge 3$.
\end{enumerate}
If, in addition, $S$ has characteristic $0$, they are also equivalent to
\begin{enumerate}[{\quad\rm(i)}]
\item[{\rm(v)}]
$\aq mSR-=0$ for some integer $m\ge 3$. \qed
\end{enumerate}
 \end{theorem}
\end{itchunk}

The general case of Conjecture II remains open. I should like to note that these
conjectures are about noetherian rings; they are false if one drops that hypothesis, see
Planas-Vilanova \cite{Pv}, and also \cite{An:ja}.

\begin{vista}
  The results in this section, and in the next, involve only the homology functors $\aq
  nSR-$.  In view of Proposition \ref{aq:fdim}, one can phrase many of them also in
  terms of the cohomology functors $\aqc nSR-$.
\end{vista}

\section{Regular homomorphisms}
In this section we turn to regular homomorphisms.  Regular local rings have been
encountered already in \ref{regularrings}. A (not-necessarily local) noetherian ring $S$
is said to be \emph{regular} if the local $S_\fq$ is regular for each prime ideal $\fq$ in
$S$.

A regular local ring is regular, because the regularity property localizes. This
last result is immediate from the characterization of regularity by Auslander,
Buchsbaum, and Serre referred to earlier, see \cite[(19.3)]{Ma}.
  
A homomorphism $\vf\col R\to S$ of noetherian rings is \emph{regular} if $S$ is flat over
$R$ and the ring $S\otimes_Rl$ is regular whenever $R\to l$ is a homomorphism essentially
of finite type and $l$ is a field.  If in addition $\vf$ is essentially of finite type,
then one says that $\vf$ is \emph{smooth}; an alternative terminology is
\emph{geometrically regular}.

\begin{example}
  Let $X$ be a finite set of variables. The inclusion $R\hra R[X]$ is smooth, whereas the
  inclusion $R\hra R[[X]]$ is regular.
\end{example}

\begin{example}
  An extension of fields $k\to l$ is regular if and only if it is separable; this is not
  too difficult to prove when $l$ is finitely generated as a field over $k$. See
  \cite[Chapter VII]{An:hca} for the argument in the general case.

  The issue with separability is well-illustrated in the following example: when $k$ is a
  field of characteristic $p$, and $a\in k$ does not have a $p$th root in $k$, the
  field extension $k\to l=k[x]/(x^p-a)$ is not geometrically regular: $l\otimes_kl$ is a
  local ring with nilpotents, and hence it is not regular.
\end{example}

\begin{remark}
  Note that the definition of a regular homomorphism has a different flavour when compared
  to that of a locally complete intersection homomorphism. The one for regularity is due
  to Grothendieck, and it is in line with his point of view that a homomorphism $\vf\col
  R\to S$ is deemed to have a certain property (regularity, complete intersection,
  Gorenstein, Cohen-Macaulay, et cetra), if the homomorphism is flat and its fibres have
  the \emph{geometric} version of the corresponding property.

  One does not define complete intersection homomorphisms this way for it would be too
  restrictive a notion; for instance, it would preclude surjective homomorphisms defined
  by regular sequences, because they are not flat.

  The content of the next exercise is that a complete intersection homomorphism in the
  sense of Grothendieck is locally complete intersection, as defined in
  \ref{lci:definition}.
\end{remark}

\begin{exercise}
\label{ci:grothendieck}
Let $\vf\col R\to S$ be a homomorphism of noetherian rings such that $S$ is flat over $R$.
Prove that $\vf$ is locally complete intersection if and only if for each prime ideal
$\fp$ in $R$, the fibre ring $S\otimes_R \kappa(\fp)$ is locally complete intersection.
\end{exercise}

The definitive criterion for regularity in terms of Andr\'e-Quillen homology is due to
Andr\'e.
  
\begin{theorem}
\label{regular:aqcriteria}
Let $\vf\col R\to S$ be a homomorphism of noetherian rings.

The following conditions are equivalent.
\begin{enumerate}[{\quad\rm(a)}]
\item
$\vf$ is regular.
\item
$\aq nSR-=0$ for each $n\geq 1$\,.
\item
$\aq 1SR-=0$\,.
\end{enumerate}
\end{theorem}

Once again, I will provide a proof only in the case where $\vf$ is essentially of finite
type: Under this hypothesis, arguing as in the proof of Theorem \ref{ci:aqcriteria}, one
may deduce Theorem \ref{regular:aqcriteria} from the following result.

\begin{theorem}
\label{regular:localcase}
Let $\vf\col (R,\fm,k)\to (S,\fn,l)$ be a local homomorphism, essentially of finite type.
The following conditions are equivalent.
\begin{enumerate}[{\quad\rm(a)}]
\item
$\vf$ is smooth.
\item
$\aq nSR-=0$ for each $n\geq 1$, and the $S$-module $\kahler SR$ is finite free.
\item
$\aq 1SRl=0$\,.
\end{enumerate}
Thus, when $\vf$ is smooth, one has that $\sctan \vf\simeq \kahler SR$.
\end{theorem}

\begin{proof}
  (a) $\implies$ (b).  The $R$-module $S$ is flat, so base change of $\vf$ along the
  composed homomorphism $R\xra{\vf}S\to l$ yields an isomorphism
\[
\aq nSRl \cong \aq n{S\otimes_Rl}ll\quad \text{for each $n$.}
\]
The composed map $l\to (S\otimes_Rl) \to l$ equals $\id^l$, so Exercise \ref{retract}
yields isomorphisms
\[
 \aq n{S\otimes_Rl}ll \cong  \aq {n+1}l{S\otimes_Rl}l \,.
\]
Since $R\to l$ is essentially of finite type, smoothness of $\vf$ implies the ring
$S'=S\otimes_Rl$ is regular. Let $\fn'$ be the maximal ideal of $S'$ such that
$S'/\fn'=l$.  Then, the local ring $S'_{\fn'}$ is regular, so Corollary \ref{regularrings}
implies the second isomorphism below
\[
\aq nl{S'}l = \aq nl{S'_{\fn'}}l = 0 \quad \text{for $n\ge 2$,}
\]
while the first one is by \ref{localization}.  Combining this with the preceding displays
yields $\aq nSRl=0$ for $n\ge 1$.  It remains to recall Lemma \ref{vanishing:local}.

Evidently, (b) $\implies$ (c).

(c) $\implies$ (a). Since $\aq 1SRl=0$, it follows from Lemma \ref{vanishing:local}
that $\sctan \vf$ is homotopy equivalent to a complex of finite free $S$-modules $L$ with
$L_1=0$.  Therefore, one has that
\[
\aqc 1SRS=\HH {-1}{\Hom SLS} = 0
\]
This is equivalent to the statement that any $R$-algebra extension of $S$ by a square-zero
ideal is split; see \ref{low degrees}.  According to a theorem of Grothendieck, this
property characterizes the smoothness of $S$ smooth over $R$; see \cite{Gr}.
\end{proof}

\begin{itchunk}{A local-global principle.}
  Let $\vf\col R\to S$ be a homomorphism of noetherian rings, $\fq$ a prime ideal in $S$,
  and set $\fp = \fq\cap R$.  One says that $\vf$ is \emph{regular at $\fq$} if $\vf_\fq$
  is flat and the $k(\fp)$-algebra $(S\otimes_R k(\fp))_\fq$ is geometrically regular.

  The exercise below is an important local-global principle for regularity. In it, the
  hypothesis that $\vf$ is essentially of finite type is insurmountable; see Remark
  \ref{cimap:completion}. There is an analogue for the complete intersection property; see
  \cite[\S5]{Av:am}.

\begin{exercise}
\label{regularity:local}
Let $\vf\col (R,\fm,k)\to (S,\fn,l)$ be local homomorphism, essentially of finite type.
Prove that if $\vf$ is regular at $\fn$, then $\vf$ is regular.
\end{exercise}
\end{itchunk}

%   In the same vein, $\vf$ is \emph{complete intersection at $\fq$} if $\vf$ is flat and
%   the local ring $(S\otimes_R k(\fp))_\fq$ is complete intersection; one does not have to
%   worry about the geometric version of the property because of Exercise
%   \ref{ci:basechange}.

Given Theorems \ref{ci:aqcriteria} and \ref{regular:aqcriteria}, it is not hard to prove
the following result, which is a crucial step in the Hochschild-Kostant-Rosenberg theorem
that calculates the Hochschild homology and cohomology of smooth algebras, see \cite{HKR},
\cite[(1.1)]{AI:tata}.

\begin{theorem}
\label{hkr}
Let $\eta\col K\to S$ be a homomorphism of noetherian rings essentially of finite
type, such that $S$ is flat as an $K$-module.

Then $\eta$ is smooth if and only if the product map $\mu^S_R\col S\otimes_KS\to S$ is
locally complete intersection.
\end{theorem}

\begin{proof}
  Set $S^e=S\otimes_KS$.  Since $\eta$ is essentially of finite type, the ring $S^e$ is
  noetherian.  Let $\psi= \eta\otimes_KS$; it is a homomorphism of rings $S\to
  S^e$, defined by $\psi(s) = 1\otimes s$ for $s\in S$.  Since $S$ is flat over
  $K$, base change yields a homotopy equivalence of complexes of $S^e$-modules:
\[
\sctan \eta\otimes_K S \simeq  \sctan{\psi}\,.
\]
The action of $S^e$ on $\sctan \eta\otimes_KS$ is given by $(s\otimes
s')(l\otimes t) = (sl \otimes s't)$.  The composition
\[
S\xra{\psi} S^e \xra{\mu}S
\]
is the identity on $S$, so Exercise \ref{retract} yields a homotopy equivalence of
$S$-modules
\[
\sctan {\mu} \simeq \susp (\sctan \psi \otimes_{S^e}S)\,.
\]
Combining the two equivalences above, one gets the  homotopy equivalence of
$S$-modules in the following diagram
\[
\sctan {\mu} \simeq \susp \big((\sctan \eta\otimes_KS)\otimes_{S^e}S\big)
\cong  \susp (\sctan \eta\otimes_SS) = \susp \sctan \eta\,.
\]
The isomorphism is justified in the exercise below.  Therefore, on the category of
$S$-modules, one has isomorphisms
\[
\aq nS{S^e}- \cong \aq {n+1}{S^e}S - \quad \text{for each $n$.}
\]
Theorems \ref{ci:aqcriteria} and \ref{regular:aqcriteria} now provide the desired
conclusion.
\end{proof}

\begin{exercise}
\label{diagonal:trick}
Let $K\to S$ be a homomorphism of rings, set $S^e=S\otimes_KS$, and let $M$ and $N$ be
$S$-modules.  As usual, $M\otimes_K N$ has a natural structure of a (right)
$S^e$-module, with $(m\otimes n).(r\otimes s) = mr\otimes sn$.  View $S$ as an
$S^e$ module via the product map $\mu\col S^e\to S$.

Prove that the natural homomorphism of $S$-modules below is bijective:
\[
(M\otimes_K N)\otimes_{S^e} S \lra M\otimes_SN\,. 
\]
Extend this result to the case when $M$ and $N$ are complexes of $S$-modules.  Caveat:
take care of the signs.
\end{exercise}

\begin{vista}
  Andr\'e-Quillen homology does not appear in the statement of Theorem \ref{hkr}.  This
  situation is typical: Andr\'e-Quillen theory provides streamlined proofs of many results
  concerning Hochschild homology, and is sometimes indispensable, see \cite{AI:inv}.
  There is a mathematical reason for this, see \cite[(8.1)]{Qu:ams}.
\end{vista}

\begin{itchunk}{\'Etale homomorphisms.}
\label{etale}
A homomorphism $\vf\col R\to S$ of noetherian rings and essentially of finite type is said
to be \emph{\'etale} if it is smooth and unramified.

\begin{exercise}
  Let $k$ be a field, and $R$ the polynomial ring $k[x_1,\dots,x_d]$.  Let $f$ be an
  element in $R$, and set $S=R/(f)$. Find necessary and sufficient conditions on $f$ for
  the homomorphism $R\to S$ to be \'etale.
\end{exercise}

\begin{exercise}
  Formulate and prove analogues of Theorems \ref{regular:aqcriteria} and \ref{hkr} for
  \'etale homomorphisms.
\end{exercise}

If you want to check whether you are on the right track, see \cite[(5.4), (5.5)]{Qu:ams}.
\end{itchunk}

\end{document}